\documentclass[UTF-8]{article}
\usepackage{graphicx}
\usepackage{authblk}  
\usepackage{bm} 
\usepackage{amsthm}
\usepackage{geometry}
\usepackage{mathrsfs}
\usepackage{algorithm}
\usepackage{amssymb}
\usepackage{algpseudocode}  
\usepackage{amsmath}
\allowdisplaybreaks[4] 
\usepackage{graphicx}
\usepackage{tikz}
\usepackage{booktabs}
\usepackage{multirow}       
\usepackage{siunitx}  
\usetikzlibrary{shapes.geometric, arrows, positioning, calc}
\usepackage{xcolor}
\numberwithin{equation}{section}
\newtheorem{definition}{Definition}
\newtheorem{assumption}{Assumption}
\newtheorem{example}{Example}
\newtheorem{lemma}{Lemma}
\newtheorem{proposition}{Proposition}
\newtheorem{theorem}{Theorem}
\newtheorem{remark}{Remark}
\newtheorem{corollary}{Corollary}

\DeclareMathOperator*{\argmax}{arg\,max}
\DeclareMathOperator*{\argmin}{arg\,min}

\title{Bayesian distributionally robust variational inequalities: regularization and quantification}

\author[1,2]{Wentao Ma\thanks{Email: \texttt{wentao.ma@connect.polyu.hk}}}
\author[1]{Zhiping Chen\thanks{Email: \texttt{zchen@mail.xjtu.edu.cn}}}
\author[2]{Xiaojun Chen\thanks{Email: \texttt{maxjchen@polyu.edu.hk}}}

\affil[1]{School of Mathematics and Statistics, Xi’an Jiaotong University, Xi’an, China}
\affil[2]{Department of Applied Mathematics, The Hong Kong Polytechnic University, Hong Kong, China}
\date{}

\begin{document}

    \maketitle
\begin{abstract}
We propose a Bayesian distributionally robust variational inequality (DRVI) framework that models the data-generating
distribution through a finite mixture family, which allows us to study the DRVI on a tractable finite-dimensional parametric ambiguity set.
To address distributional uncertainty, we construct a data-driven ambiguity set with posterior coverage guarantees via Bayesian inference. We also employ a regularization approach to ensure numerical stability.
We prove the existence of solutions to the Bayesian DRVI and the asymptotic convergence to a solution as sample size grows to infinity and the regularization parameter goes to zero. Moreover, we derive quantitative stability bounds and finite-sample guarantees under data scarcity and contamination. Numerical experiments on a distributionally robust multi-portfolio Nash equilibrium problem validate our theoretical results and demonstrate the robustness and reliability of Bayesian DRVI solutions in practice.
\\

\textbf{Keywords:} distributional robustness; Bayesian asymptotics; data-driven;  regularization method; quantitative stability
\end{abstract}

\section{Introduction}
Let $\mathcal{X} \subseteq \mathbb{R}^d$ be a nonempty closed convex set. {\color{black} For any $x\in \mathcal{X}$, let
$\mathcal{N}_{\mathcal{X}}(x):=\{v\in\mathbb{R}^d|\langle v,y-x\rangle \le 0,\ \forall y\in\mathcal{X}\}$
denote the normal cone to $\mathcal{X}$ at $x$.} The stochastic variational inequality (SVI) \cite{gurkan1999sample,robinson1996analysis} aims to find a vector \(x^* \in \mathcal{X}\) such that
\begin{equation}
		0 \in \mathbb{E}_{P^c}[\Phi(x^*, \xi)]+\mathcal{N}_\mathcal{X}(x^*), \label{SVI}
\end{equation}
where $\xi:\Omega \to \mathbb{R}^{m}$ is a random variable defined on a probability space $(\Omega,\mathcal{F},\mathbb{P})$ whose true distribution $P^c:=\mathbb{P}\circ\xi^{-1}$ is supported on $\Xi \subseteq \mathbb{R}^{m}$. The mapping $\Phi: \mathcal{X} \times \Xi \to \mathbb{R}^d$ is continuous in $x$ for any fixed $\xi\in\Xi$ and measurable and integrable in $\xi$ for any fixed $x\in\mathcal{X}$. We use $\mathbb{E}_{P^c}$ to emphasize that the expectation is with respect to $P^c$.
{\color{black}In many applications of interest, particularly in multi-agent equilibrium problems, solutions are characterized by a coupled system of best-response optimality conditions and generally cannot be reformulated as the minimization of a single objective \cite{facchinei2003finite}.
Such coupled equilibrium systems are naturally modeled within the SVI framework, with representative examples including traffic equilibria \cite{chen2017two,ordonez2010wardrop}, oligopolistic competition in the crude oil market \cite{jiang2019regularized}, and forward-contracting equilibria in electricity markets \cite{shanbhag2011complementarity}.
}

In practice, the true distribution $P^c$ is unknown.
Decision-makers therefore seek solutions that remain reliable under distributional uncertainty \cite{cutler2023stochastic}, which may stem from policy shifts, macroeconomic cycles, or market dynamics \cite{guo2023data}.
Distributionally robust optimization (DRO) addresses this challenge by replacing $P^c$ with an ambiguity set $\mathcal{P}$ designed to contain $P^c$.
This approach is well-founded in risk-sensitive decision-making, as many coherent risk measures admit dual representations closely related to DRO formulations~\cite{rockafellar2024distributional,mehta2023stochastic}.
The combination of DRO with SVI yields the distributionally robust variational inequality (DRVI) problem~\cite{jiang2024distributionally,sun2023distributionally}, which aims to find a pair $(x^*,P^*)$ such that
\begin{equation}
    \begin{aligned}
        	0 &\in \mathbb{E}_{P^*}[\Phi(x^*, \xi)]+\mathcal{N}_\mathcal{X}(x^*),\\
         P^*&\in\underset{{P}^\prime\in\mathcal{P}}{\argmax} \mathbb{E}_{{P}^\prime}[\varphi(x^*, \xi)],
    \end{aligned}\label{drvi}
\end{equation}
where $\varphi:\mathcal{X}\times\Xi\to\mathbb{R}$ is continuous in $x$ for any fixed $\xi\in\Xi$ and measurable and integrable in $\xi$ for any fixed $x\in\mathcal{X}$. Representative choices of $\varphi$ can be found in~\cite{jiang2024distributionally,sun2023distributionally}.
Common choices for ambiguity sets include those based on first- and second-order moments \cite{delage2010distributionally}, unimodality \cite{hanasusanto2015distributionally}, $\phi$-divergence~\cite{coppens2023ordered,lam2019recovering} and Wasserstein distance \cite{gao2024wasserstein,kuhn2019wasserstein}. Despite their theoretical superiority in robustness and interpretability, the construction and computational tractability of such infinite-dimensional ambiguity sets remain challenging and often lead to the curse of dimensionality and prohibitively large-scale optimization problems \cite{bertsimas2022two}.

As demonstrated in prior work (e.g., \cite{lu2013confidence,shapiro2023bayesian,wu2018bayesian}), it is natural to model $\xi$ via a parametric family of distributions $\mathcal{P}_\Theta:=\{P_{\theta}\mid\theta\in\Theta\}$, where each distribution $P_{\theta}$ is uniquely determined by a parameter $\theta$ through the probability density function (pdf) $f(\cdot|\theta)$. Specifically, for any Borel set $A\subseteq\Xi$, $P_\theta(A)=\int_A f(\xi\mid\theta)d\xi$. Here, $\Theta\subseteq\mathbb{R}^n$ is a nonempty closed parameter set. Throughout the main development, we adopt a well-specified parametric setting and assume that there exists $\theta^c\in\Theta$ such that $P^c=P_{\theta^c}$. In misspecified settings, $\theta^c$ may instead be interpreted as the parameter associated with the Kullback--Leibler (KL) projection of $P^c$ onto $\mathcal{P}_\Theta$ \cite{shapiro2023bayesian}.
Such parametric assumptions are common in SVI applications, including Gamma and Gaussian families in game theory \cite{bichler2025computing,liu2024bayesian}, uniform distributions in system control \cite{fang2007stochastic}, and Beta distributions in traffic flow management \cite{chen2012stochastic}.
However, real-world data often exhibit multimodality, skewness, and heavy-tailed behavior \cite{zhu2009worst}, which may be inadequately captured by a single parametric family. To enhance modeling flexibility while retaining a parametric structure, we work with a finite-mixture family of distributions. Theoretically, under mild regularity conditions, any continuous pdf can be approximated arbitrarily well by a finite mixture of Gaussian kernels \cite{everitt2013finite,mclachlan2000finite}. Practically, mixture models have been applied to various areas, including asset returns in portfolio selection \cite{chen2018data,zhu2014portfolio}, rewards in multi-armed bandit problems \cite{urteaga2018nonparametric}, and customer preference heterogeneity in marketing analytics \cite{rossi2003bayesian}. Thus, the finite-mixture family serves both as a structured parametric representation and as a low-dimensional approximation class.
Accordingly, we consider the finite-mixture family $P_{\theta} = \sum_{j=1}^n \theta_j Q_j$, where $\{Q_j\}_{j=1}^n$ are fixed, prespecified component distributions and $\Theta$ is the probability simplex, i.e.,  $\Theta:=\{\theta\mid e^\top\theta=1,\ \theta\in\mathbb R_+^n\}$ where $e\in\mathbb R^n$ denotes the all-ones vector.
Let $\theta^c = (\theta_1^c, \ldots, \theta_n^c)$ be the corresponding mixture weight vector and $\mathbb{B}(\theta^c, r^c):=\{\theta\mid\|\theta-\theta^c\|_\infty\leq r^c\}\subseteq\mathbb{R}^n$ be an $\ell_\infty$ norm ball centered at $\theta^c$ with radius $r^c>0$.
Motivated by parametric DRO models \cite{gupta2019near,zhu2009worst},
we propose the Bayesian DRVI formulation as
 \begin{equation}\tag{$\text{B}^c$}
    \begin{aligned}
        	0 &\in \sum_{j=1}^{n}\theta_j\mathbb{E}_{Q_j}[\Phi(x, \xi)]+\mathcal{N}_\mathcal{X}(x),\\
         \theta&\in\underset{{\theta}^\prime\in\Theta^c}{\argmax} \sum_{j=1}^{n}\theta^\prime_j\mathbb{E}_{Q_j}[\varphi(x, \xi)],
    \end{aligned}\label{bdrvi-true}
\end{equation}
where $\Theta^c:=\mathbb{B}(\theta^c,r^c)\cap\Theta$ is the true parametric ambiguity set.
{\color{black}This parametric ambiguity set is defined over mixture weights in a fixed-dimensional parameter set, which is computationally attractive when compared with ambiguity sets whose effective dimension grows with the sample size. This computational benefit comes with the need to prescribe or estimate the component distributions $Q_1,\ldots,Q_n$; see, e.g., \cite{cao2023revenue,li2019product,van2022price}. For complex unstructured data such as language, one may instead consider DRO formulations built around application-specific topic or subpopulation structure \cite{oren2019distributionally}.
}
The prespecified radius $r^c$ captures the decision-maker’s tolerance for model perturbation or risk aversion.
In particular, when $r^c=0$, $\Theta^c$ degenerates to $\{\theta^c\}$ and \eqref{bdrvi-true} reduces to the risk-neutral SVI \eqref{SVI} under the distribution $P_{\theta^c}$.

The true ambiguity set $\Theta^c$ in \eqref{bdrvi-true} is determined by $\theta^c$, which is fixed but unknown and must be estimated from limited sample data $\mathcal{S}^N=\{\xi^1,\dots,\xi^N\}$, {\color{black}where $\xi^1,\ldots,\xi^N$ are i.i.d. samples generated from the true distribution
$P_{\theta^c}$.}
Conventional approaches construct an approximation by replacing $\theta^c$ with the maximum likelihood estimate (MLE) $\hat{\theta}_N$ \cite{schottle2009robustness} and then taking the set $\mathbb{B}(\hat{\theta}_N,r^c)\cap\Theta$. However, this MLE approximation can be unstable due to finite-sample noise and model sensitivity, which may undermine robustness. To address this, we adopt a Bayesian perspective that explicitly quantifies parametric uncertainty through posterior distributions \cite{shapiro2023bayesian}.
In the Bayesian interpretation adopted here, ${\theta}^c$ is treated as a realization of a belief random variable, whose prior pdf $p({\theta})$ represents our belief before observing data \cite{wu2018bayesian}. As new data arrive, this prior is updated recursively to refine our uncertainty quantification for $\theta^c$. Formally, applying Bayes’ rule to ${\mathcal{S}^{N}}$ yields the pdf of the posterior distribution $\mathbb{P}_{\mathcal{S}^N}$ as:
\begin{equation}
	{\color{black}p({\theta}|{\mathcal{S}^{N}})=\frac{f(\mathcal{S}^N|{\theta})p({\theta})}{\int_\Theta f(\mathcal{S}^N|{\theta'})p({\theta'})d{\theta'}},}\label{new2.1old}
\end{equation}
where $f(\mathcal{S}^N|\theta):=\prod_{i=1}^Nf({\xi}^i|\theta)$.
This posterior distribution coherently integrates prior knowledge--whether from historical data, expert judgment, or physical models--and new observations, providing an evolving quantification of parametric uncertainty. In the absence of informative priors, noninformative choices such as Jeffreys’ or uniform priors are commonly used.
{\color{black}
Recent studies \cite{jaiswal2023bayesian,shapiro2025episodic,wu2018bayesian}  replace ambiguity sets with Bayesian averaging to characterize uncertainty in the model parameter through the posterior distribution.
Further, \cite{shapiro2023bayesian} constructs an ambiguity set around the corresponding parametric reference distribution using KL divergence, thereby adopting a worst-case perspective on Bayesian averaging to hedge against misspecification of the assumed parametric model.
These methodologies lead to different tradeoffs. Bayesian averaging is typically less conservative but may lack robustness,  whereas \cite{shapiro2023bayesian} provides additional protection at the cost of greater computational complexity.
Our approach can be viewed as a compromise between these two perspectives: it retains a worst-case treatment of uncertainty, but does so within the parametric family. This allows the framework to preserve tractability while integrating naturally with the DRVI structure.
}

Specifically, in the spirit of \cite{gupta2019near}, we adopt the framework of posterior-based ambiguity assessment to address parametric uncertainty in the Bayesian DRVI problem. Our goal is to construct a data-driven ambiguity set \(\hat{\Theta}_N \subseteq \Theta\) that captures parametric uncertainty as reflected by the prior knowledge and the dataset $\mathcal{S}^N$ (equivalently, by the posterior $\mathbb{P}_{\mathcal{S}^N}$). Following \cite{kuhn2019wasserstein,lam2019recovering}, we seek $\hat{\Theta}_N:=\mathbb{B}(\hat{\theta}_{N},\hat{r}_N)\cap\Theta$ that contains the true ambiguity set $\Theta^c$ with high probability. Specifically, for a prespecified posterior confidence level $\alpha\in(0,1)$, we specify nominal parameter $\hat{\theta}_N$ and radius $\hat r_N$ such that $\mathbb{P}_{\mathcal{S}^N}\left\{\Theta^c\subseteq\hat{\Theta}_N\right\}\geq1-\alpha$. To avoid overly conservative decisions in the robust formulation \cite{chan2024distributional}, $\hat r_N$ is chosen as the smallest value that satisfies the posterior coverage guarantee:
$$
\hat{r}_N=\min\left\{r\big| r\geq0,\ \mathbb{P}_{\mathcal{S}^N}\left\{\Theta^c\subseteq\mathbb{B}(\hat{\theta}_{N},r)\cap\Theta\right\}\geq1-\alpha\right\}.
$$
Based on this construction, replacing \(\Theta^c\) by its approximation \(\hat{\Theta}_N\) yields the data-driven Bayesian counterpart of \eqref{bdrvi-true}:
\begin{equation}
    \begin{aligned}
        	0 &\in \sum_{j=1}^n\theta_j\mathbb{E}_{Q_j}[\Phi(x,\xi)]+\mathcal{N}_\mathcal{X}(x),\\
         \theta&\in\underset{{\theta}^\prime\in\hat{\Theta}_N}{\argmax} \sum_{j=1}^n{\theta}^\prime_j\mathbb{E}_{Q_j}[\varphi(x,\xi)].
    \end{aligned}\label{bdrvi2}\tag{$\text{B}_N$}
\end{equation}
To address the potential numerical instability arising from the non-uniqueness of solutions to the lower-level problem in \eqref{bdrvi2}, let $\|\cdot\|$ denote the Euclidean norm and $\lambda > 0$ be a regularization parameter. We introduce the following regularized variant of~\eqref{bdrvi2}:
\begin{equation}
    \begin{aligned}
        	0 &\in \sum_{j=1}^n\theta_j\mathbb{E}_{Q_j}[\Phi(x,\xi)]+\mathcal{N}_\mathcal{X}(x),\\
         \theta&\in\underset{{\theta}^\prime\in\hat{\Theta}_N}{\argmax} \sum_{j=1}^n{\theta}^\prime_j\mathbb{E}_{Q_j}[\varphi(x,\xi)]-\lambda\|{\theta}^\prime\|^2.
    \end{aligned}\label{bdrvi-regular2}\tag{$\text{B}_N^\lambda$}
\end{equation}
Building on these formulations and the Bayesian construction, the main contributions of this paper are summarized as follows.
\begin{itemize}
\item \textbf{Modeling}.
We propose a Bayesian DRVI framework that integrates Bayesian inference with DRVI to unify robustness and flexibility in a data-driven setting.
By modeling the unknown true distribution within a finite mixture family, {\color{black}the resulting ambiguity set is defined over mixture weights and remains fixed-dimensional as sample size grows, which facilitates posterior calibration and computation.
}
The framework captures parametric uncertainty through Bayesian updating, which integrates prior knowledge with observed data. Moreover, we introduce a regularization term in the lower-level maximization problem to improve numerical stability.

\item \textbf{Analysis}.
We show that the lower-level maximization problem over the Bayesian parametric ambiguity set can be reformulated as a risk-averse formulation with Bayesian risk measures \cite{shapiro2025episodic,wu2018bayesian}.
We prove that the solutions of \eqref{bdrvi2} and \eqref{bdrvi-regular2} converge to a solution of \eqref{bdrvi-true} as $N\to\infty$ and $\lambda\downarrow0$.
Moreover, we show that the unique maximizer of the regularized lower-level problem in \eqref{bdrvi-regular2} converges to the least-norm maximizer of the unregularized problem in \eqref{bdrvi2} as $\lambda\downarrow0$.
We further derive finite-sample guarantees and quantitative stability bounds that account for both data perturbation and regularization effects, addressing practical challenges in data-driven Bayesian DRVI estimation.

\item \textbf{Applications}.
We apply the proposed Bayesian DRVI framework to a distributionally robust multi-portfolio Nash equilibrium problem, building on Bayesian parametric ambiguity sets as in \cite{zhu2014portfolio,zhu2009worst}.
Our numerical results show that (i) the method yields more accurate equilibrium solutions and lower utility error in finite samples than classical DRVI and risk-neutral SVI benchmarks, (ii) the data-driven equilibrium solutions converge to a solution of the Bayesian DRVI problem, and (iii) the method achieves favorable out-of-sample performance, balancing payoff and tail-risk mitigation. These findings underscore the practical relevance and applicability of our method in real-world decision-making scenarios.
\end{itemize}

The paper is structured as follows. Section 2 introduces a motivating Bayesian DRVI example and explores the construction of data-driven ambiguity sets from the frequentist and Bayesian perspectives. Section 3 proves the existence of solutions to the data-driven Bayesian DRVI problem and analyzes the convergence of solutions of \eqref{bdrvi-regular2} to a solution as $N\to\infty$ and $\lambda\downarrow0$. Section 4 derives quantitative stability bounds, providing finite-sample guarantees and robustness assessments under varying data regimes. Section 5 uses a distributionally robust multi-portfolio Nash equilibrium problem to show numerically the efficiency of the data-driven Bayesian DRVI with an extragradient method.

\section{Model setup and background}
We begin by illustrating how Bayesian parametric ambiguity sets can be incorporated into the DRVI framework. The motivations and theoretical foundations of DRVI have been discussed in \cite{jiang2024distributionally,sun2023distributionally}. For completeness, we provide an illustrative example below to show how Bayesian perspectives arise naturally in optimization and equilibrium problems.

\begin{example}\label{ex-final}
Finite mixture models provide a flexible framework for modeling complex distributions such as asset returns \cite{peel2000robust}. Two common interpretations of mixture models in finance are as follows: (i) if $Q_1,\ldots,Q_n$ represent asset return distributions under different market regimes from scenario analysis, then $\theta_1,\ldots,\theta_n$ are the occurrence probabilities of these regimes; (ii) if $Q_1,\ldots,Q_n$ reflect asset returns generated by different predictive models, then $\theta_1,\ldots,\theta_n$ quantify the  credibility weights of the models.
Let $\xi\sim\sum_{j=1}^n\theta^c_jQ_j$ be the random asset return vector, and let $x\in\mathcal{X}$ denote the corresponding portfolio weights. Following \cite{zhu2009worst}, the distributionally robust portfolio selection problem can be formulated as:
\begin{equation}\label{BDRO}
\min _{x \in \mathcal{X}} \max _{\theta\in\hat{\Theta}_{N}} \sum_{j=1}^n\theta_j\mathbb{E}_{Q_j}[\varphi(x,\xi)],
\end{equation}
where the loss (e.g., disutility) function $\varphi:\mathcal{X}\times\Xi\to \mathbb{R}$ is convex and continuous in \(x\). A pair $(x^*, \theta^*) \in \mathcal{X} \times \hat{\Theta}_{N}$ is a saddle point of \eqref{BDRO} if and only if the following conditions hold:
\begin{equation}\label{bi-bdro}
x^* \in \argmin _{x \in \mathcal{X}} \sum_{j=1}^n\theta^*_j\mathbb{E}_{Q_j}[\varphi(x, \xi)] \text { and } \theta^* \in \argmax _{\theta \in \hat{\Theta}_{N}} \sum_{j=1}^n\theta_j\mathbb{E}_{Q_j}[\varphi(x^*, \xi)].
\end{equation}
Assume that $\varphi$ is continuously differentiable in $x$ and the expectation $\mathbb{E}_{Q_j}[\nabla_x\varphi(x,\xi)]$ is well-defined for all $j=1,\ldots,n$, then \eqref{bi-bdro} yields the following Bayesian DRVI formulation:
$$
 0 \in \sum_{j=1}^n\theta_j^*\mathbb{E}_{Q_j}\left[\Phi\left(x^*, \xi\right)\right]+\mathcal{N}_{\mathcal{X}}\left(x^*\right) \text { and }  \theta^* \in \underset{{\theta}^\prime\in\hat{\Theta}_{N}}{\argmax } \sum_{j=1}^n{\theta}^\prime_j\mathbb{E}_{Q_{j}}\left[\varphi\left(x^*,\xi\right)\right],
$$
with $\Phi(x,\xi):=\nabla_x\varphi(x,\xi)$.
Now consider the multi-portfolio case with \(I\) accounts as in~\cite{lampariello2021equilibrium}. Let $x^i\in\mathcal{X}_i \subseteq \mathbb{R}^{d_i}$ denote the portfolio for account \(i\), and let \(x^{-i}\) represent the other \(I-1\) accounts' allocations. Let \(\varphi_i:\mathcal X_1\times\ldots\times\mathcal X_I\times\Xi\to\mathbb{R}\) be the loss function for account \(i\).  Inspired by \cite{liu2024bayesian}, a distributionally robust multi-portfolio Nash equilibrium aims to find $\left(x^{1*}, \ldots, x^{I*}\right)$ together with worst-case parameters $\theta^{i*}\in \hat{\Theta}_N$ for each account $i$ such that
\begin{equation}\label{eq-ex-7}
\begin{aligned}
& x^{i*} \in
   \argmin_{x^i \in \mathcal{X}_i}
   \sum_{j=1}^n \theta^{i*}_j  \mathbb{E}_{Q_j}\left[\varphi_i\left(x^i, x^{-i*}, \xi\right)\right], \quad i=1,\ldots,I,\\
& \theta^{i*} \in
   \argmax_{\theta^i \in \hat{\Theta}_{N}}
   \sum_{j=1}^n \theta^i_j  \mathbb{E}_{Q_j}\left[\varphi_i\left(x^{i*}, x^{-i*}, \xi\right)\right],
   \quad i=1,\ldots,I .
\end{aligned}
\end{equation}
By the same argument as for \eqref{bi-bdro},
any equilibrium $(x^{1*},\ldots,x^{I*};\theta^{1*},\ldots,\theta^{I*})$ satisfying \eqref{eq-ex-7} also satisfies the following Bayesian DRVI problem
\begin{equation}\label{eq:nash-drvi-system}
\begin{aligned}
0 &\in \sum_{j=1}^n\theta^{i*}_j\mathbb{E}_{Q_{j}}\left[\nabla_{x^i}\varphi_i\left(x^{i*}, x^{-i*}, \xi\right)\right]+\mathcal{N}_{\mathcal{X}_i}\left(x^{i*}\right), \quad i=1, \ldots, I, \\
\theta^{i*} &\in \underset{\theta^i\in\hat{\Theta}_{N}}{\argmax } \sum_{j=1}^n\theta^i_j\mathbb{E}_{Q_{j}}\left[\varphi_i\left(x^{i*}, x^{-i*}, \xi\right)\right], \quad i=1, \ldots, I,
\end{aligned}
\end{equation}
where $\Phi_i(x^1,\ldots,x^I,\xi):=\nabla_{x^i}\varphi_i(x^1,\ldots,x^I,\xi)$.
\end{example}

This example illustrates the versatility of the Bayesian DRVI framework in modeling complex decision problems, particularly those arising in finance.
{\color{black}
While the single-agent model \eqref{BDRO} can be treated as a standard min--max program,
the multi-portfolio model~\eqref{eq-ex-7} is inherently an equilibrium system: it is defined by coupled optimality conditions across agents,
and in general does not admit a reformulation as the minimization of a single objective.
This motivates adopting Bayesian DRVI as a natural representation for robust multi-agent equilibrium systems.}
Below, we detail the construction of $\hat{\Theta}_N$.

\subsection{Bayesian parametric ambiguity sets for distributions}

\textcolor{black}{
We work with a Bayesian parametric ambiguity set defined over the mixture weights of a prescribed finite-mixture family. As a result, the ambiguity set is constructed in a parameter set whose dimension remains fixed as sample size grows, in contrast to data-driven ambiguity sets supported on the observed samples \cite{lam2019recovering,kuhn2019wasserstein}, whose effective dimension typically scales with the number of observations. This representation is computationally attractive, but it comes with a stronger modeling commitment: the component distributions $Q_1,\ldots,Q_n$ must be prescribed from structural knowledge or estimated from data.
Hence, the proposed framework is most natural in settings where such mixture components have meaningful structural interpretations or can be reliably estimated.
In practice, the components may reflect regimes, latent segments, or distinct demand-generation mechanisms \cite{cao2023revenue,li2019product,van2022price}, or they may be learned using standard mixture-estimation procedures such as the expectation--maximization (EM) algorithm \cite{tseng2004analysis,redner1984mixture}.}
To construct practical ambiguity sets that contain the true parametric set $\Theta^c$ with high probability, we first introduce two fundamental concepts of confidence regions with coverage guarantees in the DRO context.

\begin{definition}[Frequentist guarantee]\label{def-1}
A set $\mathcal{P}_N\subseteq\mathcal{P}$ is called a level-\(\alpha\) frequentist confidence region if  $\mathbb{P}_{\mathcal{S}^N | P^c}\left\{P^c\in\mathcal{P}_{N}\right\} \geq 1-\alpha$, where $\mathbb{P}_{\mathcal{S}^N | P^c}:=(P^c)^N$ denotes the probability measure induced by $N$ i.i.d. observations from $P^c$.
\end{definition}

Definition \ref{def-1} states that a frequentist confidence region $\mathcal{P}_N$ contains the true distribution $P^c$ with probability at least $1-\alpha$ under repeated sampling.  {\color{black}Classical ambiguity sets satisfying such guarantees include constructions based on moments \cite{delage2010distributionally}, $\phi$-divergences \cite{ben2013robust,coppens2023ordered,shapiro2017distributionally}, hypothesis tests \cite{bertsimas2018robust}, and Wasserstein/$\ell_1$/$\ell_\infty$ metrics \cite{jiang2018risk,kuhn2019wasserstein}. In our setting, distributional uncertainty is modeled within a prescribed finite-mixture family through a Bayesian parametric ambiguity set over the mixture weights. This yields a fixed-dimensional parametric representation, while the uncertainty quantification is expressed in posterior terms. We therefore introduce below a Bayesian analogue of Definition \ref{def-1} in the parametric setting.}

\begin{definition}[Posterior guarantee]\label{def-2}
A set $\mathcal{P}_{\Theta_N}:=\{P_\theta\mid\theta\in\Theta_N\}\subseteq\mathcal{P}_\Theta$ is called a level-$\alpha$ Bayesian confidence region if
$\mathbb{P}_{\mathcal{S}^N}\left\{P_{\theta^c}\in\mathcal{P}_{\Theta_N}\right\} \geq 1-\alpha$.
\end{definition}

In contrast to Definition~\ref{def-1}, the probability here is taken with respect to the posterior distribution $\mathbb{P}_{\mathcal{S}^N}$ induced by the observed data.
The Bayesian confidence region $\mathcal{P}_{\Theta_N}$ in Definition~\ref{def-2} corresponds one-to-one to a parameter set $\Theta_N\subseteq\Theta$, which is commonly referred to as a level-$\alpha$ Bayesian confidence set in \cite{lehmann2005testing}, satisfying
$\mathbb{P}_{\mathcal{S}^N}\left\{{\theta^c}\in{\Theta_N}\right\} \geq 1-\alpha$.
For notational simplicity, we henceforth work directly with the parameter set $\Theta_N$.
Moreover, any $\Theta_N$ that satisfies the posterior guarantee at level $\alpha$ automatically satisfies the Bayesian feasible guarantee in \cite[Definition 2]{gupta2019near} at the same confidence level. For a detailed discussion on the trade-offs between frequentist and Bayesian methods, we refer the reader to  \cite{efron2021computer,gupta2019near}.

As multiple sets may satisfy Definition \ref{def-2}, we aim to construct a computationally tractable minimal Bayesian confidence region $\Theta_N$, thereby avoiding excessive conservatism \cite{chan2024distributional}.
{\color{black}To obtain an explicit construction, we invoke a normal approximation of the posterior under standard Bernstein--von Mises regularity conditions \cite{gupta2019near}.
\begin{assumption}\label{ass-bvm}
Assume (i) the true parameter $\theta^c$ lies in the relative interior of $\Theta$, (ii) the log-likelihood $\log f(\mathcal S^N|\theta)$ is locally twice continuously differentiable around $\theta^c$, (iii) the Fisher information matrix ${\cal I}(\theta^c)$ is nonsingular, and
(iv) the prior pdf $p(\theta)$ is continuous and strictly positive in a neighborhood of $\theta^c$.
\end{assumption}
}
Under Assumption~\ref{ass-bvm}, the posterior $\mathbb{P}_{\mathcal{S}^N}$ can be approximated by a normal distribution $\mathcal{N}(\hat{{\theta}}_N,\Sigma_N)$ by the Bernstein--von Mises theorem \cite{van2000asymptotic}, with
$\Sigma_N =\frac{1}{N}{\cal I}(\hat\theta_N)^{-1}$ and ${\cal I}(\hat\theta_N)=-\frac{1}{N}\nabla^2_\theta\log f(\mathcal S^N|\theta)\Big|_{\theta=\hat\theta_N}$. This approximation is widely employed in the literature; see, e.g., \cite{chen2018data,gupta2019near,zhu2014portfolio}.
Motivated by this and applying a Bonferroni correction for simultaneous marginal coverage, we define the confidence radius
$
\hat\delta_N
=
\max_{1\le j\le n}
\sqrt{\Sigma_{N,jj}}
z_{1-\frac{\alpha}{2n}},
$
where $\Sigma_{N,jj}$ is the $j$-th diagonal element of $\Sigma_N$ and
$z_{1-\frac{\alpha}{2n}}$ is the $(1-\frac{\alpha}{2n})$-quantile of the standard normal distribution.
This choice leads to the explicit construction
\begin{equation}\label{eq-confidence-region}
    \Theta_N:=\mathbb{B}( \hat{\theta}_N,\hat{\delta}_N)\cap\Theta=\bigl\{\theta\mid\|\theta-\hat\theta_N\|_\infty\le\hat\delta_N,\ \theta\in\Theta\bigr\}
\end{equation}
as a computationally convenient approximation of the minimal level-$\alpha$ Bayesian confidence region.

{\color{black}To quantify the convergence of the Bayesian confidence region $\Theta_N$, the data-driven ambiguity set $\hat\Theta_N$, and later the corresponding solution sets of \eqref{bdrvi2} and \eqref{bdrvi-regular2},} we employ the directed set deviation $\mathbb{D}$ induced by the $\ell_\infty$ norm. Specifically, for a metric space $(\mathbb{X}, \|\cdot\|_\infty)$ and any two subsets $S_1, S_2 \subseteq \mathbb{X}$, the discrepancy measure is defined as:
$$\mathbb{D}\left(S_1, S_2\right)=\sup _{s_1 \in S_1} \mathbf{d}\left(s_1, S_2\right)=\sup _{s_1 \in S_1} \inf _{s_2 \in S_2} \|s_1- s_2\|_\infty.$$
Further, the Hausdorff distance is defined as  $ \mathbb{H}\left(S_1, S_2 \right)=\max \left\{\mathbb{D}\left(S_1, S_2\right), \mathbb{D}\left(S_2, S_1\right)\right\}$.
{\color{black}
We next impose mild conditions to establish almost-sure consistency of the Bayesian confidence region constructed above, which will be used later for convergence analysis.
Following the setting in \cite{wu2018bayesian}, the observed data $\mathcal S^N=\{\xi^1,\ldots,\xi^N\}$ are generated i.i.d. from the true distribution
$P_{\theta^c}$, and the infinite data stream $(\xi^1,\xi^2,\ldots)$ induces the product probability measure
$(P_{\theta^c})^{\mathbb N}$ on the sample-path space $\Xi^{\mathbb N}$.
Unless otherwise stated, every ``almost surely'' (a.s.) convergence in this paper is with respect to
$(P_{\theta^c})^{\mathbb N}$.
}

\begin{assumption}\label{ass3.1} For all $\theta\in\Theta$, (i) there exists an integrable function $g$ such that $|\log f(\xi|\theta)|\le g(\xi)$ with $\mathbb{E}_{P_{\theta^c}}[g(\xi)]<\infty$, and (ii) there exist constants $0<c_1 < c_2$ such that the prior pdf $p(\theta)$ satisfies $c_1 \leq p({\theta}) \leq c_2$.
\end{assumption}
Under the above conditions, we can establish the convergence of the Bayesian confidence region $\Theta_N$ defined in \eqref{eq-confidence-region} as $N \to \infty$.
\begin{lemma}\label{lem-con-ambiguity}
Under Assumptions \ref{ass-bvm}-\ref{ass3.1}, $\mathbb{H}({\Theta}_N,\{\theta^c\})$ converges to 0 almost surely as $N\to\infty$.
\end{lemma}
This follows from standard results (e.g., \cite[Lemma 3.1]{shapiro2023bayesian}) and is referred to as Bayesian consistency in \cite[Assumption 3.2]{shapiro2025episodic}. {\color{black} In particular, $\|\hat\theta_N-\theta^c\|_\infty\to 0$ and $\hat\delta_N\to 0$ almost surely as $N\to\infty$.}
Based on the Bayesian confidence region $\Theta_N$ in \eqref{eq-confidence-region} and Lemma \ref{lem-con-ambiguity}, we now construct the Bayesian parametric ambiguity set $\hat{\Theta}_N$ with the minimal (i.e., smallest) radius that still contains $\Theta^c$ with posterior probability at least $1-\alpha$.

\begin{theorem}\label{thm-new}
{\color{black}Suppose Assumption~\ref{ass-bvm} holds.} Let $\hat r_N := r^c+\hat\delta_N$ and $\hat{\Theta}_N := \mathbb{B}(\hat{\theta}_N,\hat r_N)\cap\Theta$.
Then
$$\mathbb{P}_{\mathcal{S}^N}\left\{\Theta^c\subseteq\hat{\Theta}_N\right\} \geq 1-\alpha \text{ and }
\hat{r}_N=\min\left\{r\big| r\geq0,\ \mathbb{P}_{\mathcal{S}^N}\left\{\Theta^c\subseteq\mathbb{B}(\hat{\theta}_{N},r)\cap\Theta\right\}\geq1-\alpha\right\}.
$$
Furthermore, under Assumptions~\ref{ass-bvm}-\ref{ass3.1}, $\mathbb{H}(\hat{\Theta}_N,\Theta^c)$ converges to 0 almost surely as $N\to\infty$.
\end{theorem}
 \textit{Proof}.
    Recall that the ambiguity set is defined using the $\ell_\infty$ norm. Suppose that $\theta^c\in\Theta_N$, then for any $\theta\in\Theta^c$, we have
    $$
    \|\theta-\hat{\theta}_N\|_\infty\leq\|\theta-\theta^c\|_\infty+\|\theta^c-\hat{\theta}_N\|_\infty\leq r^c+\hat{\delta}_N.
    $$
    Thus, whenever $\theta^c \in \Theta_N$, the true ambiguity set $\Theta^c$ is contained in $\hat\Theta_N$.
    In terms of probability, this implies
$
\mathbb{P}_{\mathcal{S}^N}\left\{\Theta^c\subseteq\hat{\Theta}_N\right\} \geq \mathbb{P}_{\mathcal{S}^N}\left\{\theta^c\in\Theta_N\right\}\geq1-\alpha.
$
To verify the minimality of $\hat{r}_N$, assume that there exists some $\delta'_N<\hat{\delta}_N$ such that
    \[
    \mathbb{P}_{\mathcal{S}^N}\left\{\Theta^c\subseteq \mathbb{B}(\hat{\theta}_N, \delta'_N+r^c)\cap\Theta\right\} \geq 1-\alpha.
    \]
 For $\Theta^c\subseteq \mathbb{B}(\hat{\theta}_N, \delta'_N+r^c)\cap\Theta$ to hold, we must have $
    \theta^c \in \mathbb{B}(\hat{\theta}_N, \delta'_N)\cap\Theta$.
 This implies
    \[
    \mathbb{P}_{\mathcal{S}^N}\left\{\theta^c \in \mathbb{B}(\hat{\theta}_N, \delta'_N)\cap\Theta\right\} \geq \mathbb{P}_{\mathcal{S}^N}\left\{\Theta^c\subseteq \mathbb{B}(\hat{\theta}_N, \delta'_N+r^c)\cap\Theta\right\} \geq 1-\alpha,
    \]
    which contradicts the minimality of $\hat{\delta}_N$. Finally, under Assumption~\ref{ass3.1},  the convergence of $\hat{\Theta}_N$ to $\Theta^c$ follows directly from Lemma \ref{lem-con-ambiguity}.
$\hfill\blacksquare$

Theorem~\ref{thm-new} shows that, under the parametric mixture model and the Bernstein--von Mises theorem \cite{van2000asymptotic} with Bonferroni correction, the infinite-dimensional distributional ambiguity set over distributions {\color{black}can be represented by a finite-dimensional ambiguity set $\hat{\Theta}_N=\mathbb{B}(\hat\theta_N, \hat{r}_N)\cap\Theta$ over mixture weights} with a posterior coverage guarantee.

\subsection{Model component estimation and comparison with classical ambiguity sets}
{\color{black}
The construction in Section~2.1 assumes that the component distributions
$Q_1,\ldots,Q_n$ are fixed. In practice, however, these components may be estimated from data rather than specified a priori.
To study this practically relevant setting, we consider a data-driven experiment in which Gaussian components are learned by EM while the number of mixture components is varied explicitly.
The goal is to understand how the proposed Bayesian ambiguity set behaves when the component distributions are estimated from data and the number of mixture components varies.

We consider a continuous true distribution $P^c$ given by a three-component mixture of Student-$t$ distributions with unequal locations, scales, and degrees of freedom, namely
\[
P^c
=
0.55t_4(-2,0.8^2)
+
0.30t_6(1.2,0.5^2)
+
0.15t_3(4,0.7^2),
\]
where
$t_\nu(\mu,\sigma^2)$ denotes the Student-t distribution with $\nu$ degrees of freedom, location
$\mu$, and scale parameter $\sigma$.
This choice yields a multimodal, asymmetric, and heavy-tailed distribution, thereby providing a misspecified test case for Gaussian-mixture modeling.
For each prescribed number of components $n$, we fit Gaussian components $Q_1,\ldots,Q_n$ to the sample $\mathcal S^N$ by EM and then construct the Bayesian ambiguity set over the mixture weights as in Section~2.1.

We evaluate the resulting construction from two perspectives.
First, we measure the nominal approximation quality of the fitted Gaussian-mixture model through the integrated  cumulative distribution function (CDF) error
\[
\int_{\mathbb R}\bigl|\hat F_n(\xi)-F^c(\xi)\bigr|d\xi,
\]
where $\hat F_n$ denotes the CDF induced by the center of the fitted $n$-component Bayesian ambiguity set, and
$F^c$ is the CDF of the true data-generating distribution $P^c$.
Second, we quantify conservatism through the average width of the induced ambiguity band, computed from the pointwise lower and upper CDF envelopes associated with the ambiguity set.
Table~\ref{tab:misspecification-main} reports the Monte Carlo means over $100$ replications.

\begin{table}[ht]
\centering
\caption{\textcolor{black}{Integrated CDF error and average ambiguity-band width under different sample sizes $N$ and numbers of Gaussian mixture components $n$.}}
\label{tab:misspecification-main}
\begin{tabular}{cc|cccccccc}
\toprule     
& & \multicolumn{4}{c}{\textcolor{black}{Integrated CDF error}} & \multicolumn{4}{c}{\textcolor{black}{Average ambiguity-band width}} \\
\cmidrule(lr){3-6} \cmidrule(lr){7-10}
\textcolor{black}{$N$} &  & \textcolor{black}{$n=2$} & \textcolor{black}{$n=5$} & \textcolor{black}{$n=8$} & \textcolor{black}{$n=10$} & \textcolor{black}{$n=2$} & \textcolor{black}{$n=5$} & \textcolor{black}{$n=8$} & \textcolor{black}{$n=10$} \\ 
\midrule
\textcolor{black}{50}   & & \textcolor{black}{0.4214} & \textcolor{black}{0.4162} & \textcolor{black}{0.4144} & \textcolor{black}{0.4087} & \textcolor{black}{0.0757} & \textcolor{black}{0.2298} & \textcolor{black}{0.3231} & \textcolor{black}{0.3582} \\
\textcolor{black}{100}  & & \textcolor{black}{0.3252} & \textcolor{black}{0.2875} & \textcolor{black}{0.2873} & \textcolor{black}{0.2831} & \textcolor{black}{0.0540} & \textcolor{black}{0.1666} & \textcolor{black}{0.2425} & \textcolor{black}{0.2763} \\
\textcolor{black}{200}  & & \textcolor{black}{0.2993} & \textcolor{black}{0.2397} & \textcolor{black}{0.2393} & \textcolor{black}{0.2390} & \textcolor{black}{0.0353} & \textcolor{black}{0.1247} & \textcolor{black}{0.1882} & \textcolor{black}{0.2248} \\
\textcolor{black}{500}  & & \textcolor{black}{0.2194} & \textcolor{black}{0.1292} & \textcolor{black}{0.1290} & \textcolor{black}{0.1285} & \textcolor{black}{0.0219} & \textcolor{black}{0.0793} & \textcolor{black}{0.1271} & \textcolor{black}{0.1558} \\ \bottomrule
\end{tabular}
\end{table}

Table~\ref{tab:misspecification-main} shows that, for each fixed sample size $N$, increasing the number $n$ of Gaussian mixture components reduces the approximation error and enlarges the ambiguity-band width.
At the same time, for each fixed $n$, increasing the number of sample size $N$ reduces both the approximation error and the ambiguity-band width.
These results suggest that moderate mixture orders provide a favorable balance between approximation quality and ambiguity-band width in this data-driven setting.

To place the above data-driven construction in context, we also report a supplementary comparison
with classical ambiguity sets.
For this benchmark comparison, we select the number of Gaussian mixture components $n$ by the Bayesian information criterion (BIC) \cite{schwarz1978estimating} and estimate the corresponding Gaussian mixture model from the sample by the EM algorithm.}
As classical benchmarks, given samples $\mathcal{S}^N:=\{\xi^1,\ldots,\xi^N\}$ from $P^c$, the $\phi$-divergence ambiguity set (e.g., modified $\chi^2$ divergence) in \cite{ben2013robust,lam2019recovering} and the $\ell_1$ norm ambiguity set in \cite{farokhi2023distributionally,jiang2018risk,wasserman2006all} that satisfy the frequentist guarantee at level $\alpha$ can be constructed as follows:
\begin{equation*}
  \text{$\ell_1$ norm: }
  \mathcal{P}^{\ell_1}_{N}=\Biggl\{P \bigg|
  p=(p_1,\ldots,p_N)\in \Delta^N,\
  \sum_{i=1}^{N}|p_i-\hat{p}_{N,i}|\leq \sqrt{\frac{2}{N}\log\left(\frac{2}{\alpha}\right)}\Biggr\},
\end{equation*}
\begin{equation*}
  \text{Modified }\chi^2\text{ divergence: }
  \mathcal{P}^{\chi^2}_{N}=\Biggl\{P \bigg|
  p=(p_1,\ldots,p_N)\in \Delta^N,\
  \sum_{i=1}^{N}\frac{(p_i-\hat{p}_{N,i})^2}{\hat{p}_{N,i}}\leq \frac{\chi^2_{N-1,1-\alpha}}{N}\Biggr\}.
\end{equation*}
Here, $\Delta^N:=\{p\mid \sum_{i=1}^Np_i=1,\ p\in\mathbb{R}^N_+\}$ is the $N$-dimensional probability simplex and $\hat{p}_{N}:=(1/N,\ldots,1/N)$ denotes the probability vector of the empirical distribution based on $\mathcal{S}^N$.
Each distribution $P$ supported on $\mathcal{S}^N$ is uniquely identified with its probability vector
$p=(p_1,\ldots,p_N)\in\Delta^N$, where $p_i:=P(\xi=\xi^i)$ represents the probability mass at sample $\xi^i$.
The term $\chi^2_{N-1,1-\alpha}$ denotes the $(1-\alpha)$-quantile of the $\chi^2_{N-1}$ distribution with $N-1$ degrees of freedom.

\begin{figure}[h]
    \centering
    \includegraphics[width=\textwidth]{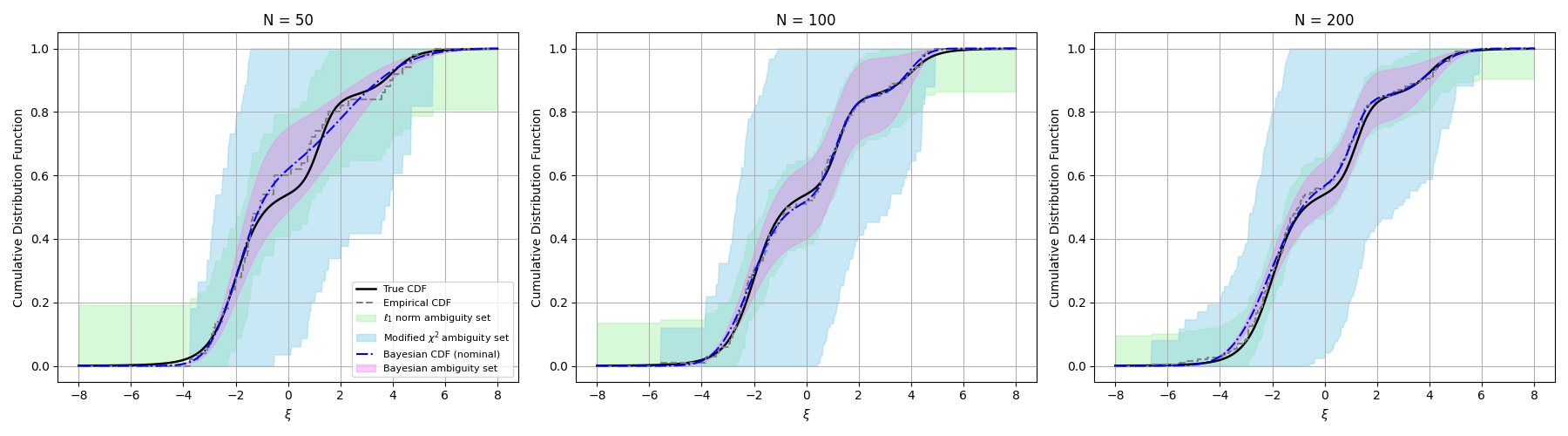}
    \caption{CDF envelopes induced by different ambiguity sets at $\alpha=0.05$.}
    \label{fig:benchmark-cdf}
\end{figure}

To visualize the geometry induced by these constructions, we plot the pointwise lower and upper
envelopes of the CDF associated with the fitted Bayesian ambiguity set
$\hat\Theta_N$ and with the two benchmark sets at level $\alpha=0.05$ for
sample sizes $N\in\{50,100,200\}$; see Figure~\ref{fig:benchmark-cdf}.

{\color{black}The CDF envelopes in Figure~\ref{fig:benchmark-cdf} suggest 
that under the fitted Gaussian-mixture model, the nominal Bayesian CDF provides
a smooth approximation to the true data-generating distribution and the ambiguity band induced by the Bayesian mixture model is typically
narrower than those induced by the $\ell_1$ and modified $\chi^2$ constructions.

}

\section{Existence and convergence analysis of solutions}
In this section, we first prove the existence of solutions to (\ref{bdrvi2}) and (\ref{bdrvi-regular2}) under mild conditions. Next, we show that as $N\to\infty$ and $\lambda\downarrow0$, the solutions of the data-driven problems (\ref{bdrvi2}) and (\ref{bdrvi-regular2}) converge almost surely to a solution of problem \eqref{bdrvi-true}.

\subsection{Existence of solutions for \eqref{bdrvi2} and \eqref{bdrvi-regular2}}
\textcolor{black}{We begin by verifying that the Bayesian DRVI problems are well posed.
Specifically, for any fixed $\hat{\theta}_N$ and $\hat{r}_N$, we show that the $x$-solution sets of the Bayesian DRVI problem (\ref{bdrvi2}) and its regularized counterpart (\ref{bdrvi-regular2}) are nonempty and compact. This provides a basic regularity framework for the convergence and stability analysis developed below.}
To this end, we first introduce the necessary notation. For any fixed $x \in \mathcal{X}$, denote the optimal value and the set of maximizers (in $\theta$) of the lower-level maximization problem in \eqref{bdrvi2}, as well as those associated with the regularized lower-level maximization problem in \eqref{bdrvi-regular2}, by
$$
\nu(x,  \hat{\theta}_N,\hat{r}_N):=\max _{\theta \in \hat{\Theta}_N} \sum_{j=1}^n\theta_j\mathbb{E}_{Q_j}[\varphi(x,\xi)],\ \nu^\lambda(x,  \hat{\theta}_N,\hat{r}_N):=\max _{\theta \in \hat{\Theta}_N} \sum_{j=1}^n\theta_j\mathbb{E}_{Q_j}[\varphi(x,\xi)]-\lambda\|\theta\|^2, $$
$$\vartheta(x,  \hat{\theta}_N,\hat{r}_N):=\underset{\theta\in\hat{\Theta}_N}{\argmax } \sum_{j=1}^n\theta_j\mathbb{E}_{Q_j}[\varphi(x,\xi)],\ \vartheta^\lambda(x,  \hat{\theta}_N,\hat{r}_N):=\underset{\theta\in\hat{\Theta}_N}{\argmax } \sum_{j=1}^n\theta_j\mathbb{E}_{Q_j}[\varphi(x,\xi)]-\lambda\|\theta\|^2,
$$
respectively.
Moreover, let ${\mathfrak{X}(\theta^c,r^c)}$, $\mathfrak{X}( \hat{\theta}_N,\hat{r}_N)$ and $\mathfrak{X}^\lambda( \hat{\theta}_N,\hat{r}_N)$ denote the $x$-solution sets of Bayesian DRVI \eqref{bdrvi-true}, the data-driven Bayesian DRVI \eqref{bdrvi2} and the data-driven regularized Bayesian DRVI \eqref{bdrvi-regular2}, respectively.

\begin{proposition}\label{prop-exist}
    Suppose that \(\mathbb{E}_{Q_j} [\Phi(\cdot, \xi)] \) and \(\mathbb{E}_{Q_j} [\varphi(\cdot, \xi)]\) are continuous for $j\in\{1,\ldots,n\}$ and the set $\mathcal{X}$ is bounded. Then $\mathfrak{X}( \hat{\theta}_N,\hat{r}_N)$ and $\mathfrak{X}^\lambda( \hat{\theta}_N,\hat{r}_N)$ are nonempty and compact.
\end{proposition}

\textit{Proof}.   The result for $\mathfrak{X}^\lambda( \hat{\theta}_N,\hat{r}_N)$ follows immediately when $\lambda>0$, together with the continuity of $\mathbb{E}_{Q_j}[\Phi(x, \xi)] $ and the compactness of $\mathcal{X}$.
We only need to focus on the case $\lambda = 0$.
For any fixed $x\in\mathcal{X}$, the maximization problem
$ \theta\in\underset{{\theta}^\prime\in\hat{\Theta}_N}{\argmax} \sum_{j=1}^n{\theta}^\prime_j\mathbb{E}_{Q_j}[\varphi(x,\xi)]$
can be reformulated as
$0 \in-{\varphi}^n(x)+\mathcal{N}_{\hat{\Theta}_N}(\theta)$, where
$
{\varphi}^n(x):=\left(\mathbb{E}_{Q_1}[\varphi(x,\xi)],\ldots,\mathbb{E}_{Q_n}[\varphi(x,\xi)]\right)^{\top} \in \mathbb{R}^{n}
$
and $\mathcal{N}_{\hat{\Theta}_N}(\theta)$ is the normal cone to the compact set $\hat{\Theta}_N$ at $\theta$. Therefore, \eqref{bdrvi2} can be written as the following VI problem:
$$
\begin{aligned}
0 \in \sum_{j=1}^n \theta_j \mathbb{E}_{Q_j}[\Phi(x,\xi)]+\mathcal{N}_{\mathcal{X}}(x) \text{ and } 0 \in-{\varphi}^n(x)+\mathcal{N}_{\hat{\Theta}_N}(\theta).
\end{aligned}
$$
Since $\mathcal{X}\times\hat{\Theta}_N$ is nonempty convex and compact and $\varphi^n$ is continuous, the desired conclusion follows directly from Corollary 2.2.5 in \cite{facchinei2003finite}.
$\hfill\blacksquare$

\subsection{Convergence analysis}
In this subsection, we prove the
convergence of the $x$-solution sets $\mathfrak{X}^\lambda( \hat{\theta}_N,\hat{r}_N)$ of problem \eqref{bdrvi-regular2} to
$\mathfrak{X}(\theta^c,r^c)$ of \eqref{bdrvi-true} as $N\to\infty$ and $\lambda\downarrow0$. The convergence analysis proceeds in two steps: (i) the convergence of
$\mathfrak{X}^\lambda( \hat{\theta}_N,\hat{r}_N)$ to $\mathfrak{X}( \hat{\theta}_N,\hat{r}_N)$ as $\lambda\downarrow0$, and (ii) the convergence of $\mathfrak{X}( \hat{\theta}_N,\hat{r}_N)$ to $\mathfrak{X}(\theta^c,r^c)$ as $N\to\infty$. Due to the continuity of  the functions and the compactness of $\mathcal{X}$, there exists \(M>0\) such that \( \sup_{x\in\mathcal{X},j\in\{1,\ldots,n\}} \{|\mathbb{E}_{Q_j}[\varphi(x,\xi)]|,|\mathbb{E}_{Q_j}[\Phi_1(x,\xi)]|,\ldots,|\mathbb{E}_{Q_j}[\Phi_d(x,\xi)]|\}\leq M\). For each $j=1,\dots,n$, we define the component-wise lower and upper bounds for the mixture weights as {\color{black}${\hat{\theta}}^l_{N,j}:=\max\{0,\ \hat{\theta}_{N,j}-\hat{r}_N\}$ and ${\hat{\theta}}^u_{N,j}:=\min\{1,\ \hat{\theta}_{N,j}+\hat{r}_N\}.$
}
The following result provides an alternative representation of the lower-level objective function, highlighting the risk-averse nature of our Bayesian DRVI framework.

{\color{black}
\begin{theorem}\label{thm-reformulation}
Let \(J\) be a categorical random variable taking values in \(\{1,\dots,n\}\).
Let $P^l$ and $P^{u-l}$ be the probability distributions of $J$ with
\[
P^l(J= j)=\frac{{\hat{\theta}}^l_{N,j}}{\sum_{i=1}^n{\hat{\theta}}^l_{N,i}}
\quad\text{and}\quad
P^{u-l}(J= {j})=\frac{{\hat{\theta}}^u_{N,j}-{\hat{\theta}}^l_{N,j}}{\sum_{i=1}^n(\hat{\theta}^u_{N,i}-{\hat{\theta}}^l_{N,i})}.
\]
Define
$\beta_N
:=
1-\frac{1-\sum_{j=1}^n{\hat{\theta}}^l_{N,j}}{\sum_{i=1}^n(\hat{\theta}^u_{N,i}-{\hat{\theta}}^l_{N,i})}$.
Then, for any fixed $x\in\mathcal{X}$, we have
\[
\nu(x,  \hat{\theta}_N,\hat{r}_N)
=
\Bigl(\sum_{j=1}^n{\hat{\theta}}^l_{N,j}\Bigr)\mathbb{E}_{P^l}\mathbb{E}_{Q_J}[\varphi(x,\xi)]
+\Bigl(1-\sum_{j=1}^n{\hat{\theta}}^l_{N,j}\Bigr)
\mathrm{CVaR}^{\beta_N}_{P^{u-l}}\mathbb{E}_{Q_J}[\varphi(x,\xi)],\]
where $\mathrm{CVaR}^{\beta_N}_{P^{u-l}}\mathbb{E}_{Q_J}[\varphi(x,\xi)] $ represents the conditional value-at-risk (CVaR) of $\mathbb{E}_{Q_J}[\varphi(x,\xi)]$ with respect to $P^{u-l}$ at risk level $\beta_N$.
\end{theorem}
}

 \textit{Proof}.
Without loss of generality, for a given $x$, we assume that the values $\varphi_j:=\mathbb{E}_{Q_j}[\varphi(x,\xi)]$ are sorted in a non-decreasing order:
$\varphi_1\leq\ldots\leq \varphi_n$.
Recall that $\nu(x,\hat{\theta}_N,\hat{r}_N)$ equals the optimal value of the following optimization problem:	\begin{equation}
		\begin{split}
			\sup_{\theta\in\mathbb{R}^n}&\sum_{j=1}^n\theta_{j}\varphi_j\\
			\text{s.t.}\ &\hat{\theta}_{N,j}^l\leq \theta_{j}\leq \hat{\theta}_{N,j}^u,\ j=1,\ldots,n,\ e^\top\theta=1.
		\end{split}\label{3.4}
	\end{equation}
Introducing Lagrange multipliers
\(d_u^{j}\ge0\) for the upper bounds,
\(d_l^{j}\ge0\) for the lower bounds,
and \(d_0\in\mathbb R\) for the equality constraint, the Lagrange dual problem of (\ref{3.4}) can be reformulated as
	\begin{align*}	    			&\inf_{d_u,d_l\geq0,d_0}\sup_{\theta\in\mathbb{R}^n}\sum_{j=1}^n\theta_{j}\varphi_j+\sum_{j=1}^nd_u^{j}(\hat{\theta}_{N,j}^u-\theta_{j})+\sum_{j=1}^nd_l^{j}(\theta_{j}-\hat{\theta}_{N,j}^l)+d_0\left(1-\sum_{j=1}^n\theta_{j}\right)\\
			=& \inf_{d_u,d_l\geq0,d_0}d_0+\sum_{j=1}^n(d_u^{j}\hat{\theta}_{N,j}^u-d_l^{j}\hat{\theta}_{N,j}^l)+\sum_{j=1}^n\sup_{{\theta}_{j}\in\mathbb{R}}(\varphi_j+d_l^{j}-d_u^{j}-d_0)\theta_{j}\\
			=&\inf_{d_0}d_0+	\sum_{j=1}^n\inf_{\substack{d_u^{j},d_l^{j}\geq0\\ \varphi_j+d_l^{j}-d_u^{j}-d_0=0}}\left\{d_u^{j}\hat{\theta}_{N,j}^u-d_l^{j}\hat{\theta}_{N,j}^l\right\}	\\
			=&\inf_{d_0}d_0+	\sum_{j=1}^{j^*}(\varphi_j-d_0)\hat{\theta}_{N,{j}}^l+\sum_{j=j^*+1}^n(\varphi_j-d_0)\hat{\theta}_{N,{j}}^u\\
			=&\inf_{d_0}d_0+	\sum_{j=1}^{n}\hat{\theta}_{N,j}^l(\varphi_j-d_0)+\sum_{j=j^*+1}^n(\hat{\theta}_{N,{j}}^u-\hat{\theta}_{N,{j}}^l)(\varphi_j-d_0)\\			=&\sum_{j=1}^n\hat{\theta}_{N,j}^l\varphi_j+\inf_{d_0}\left\{(1-\sum_{j=1}^n{\hat{\theta}}^l_{N,j})d_0+\sum_{j=1}^n(\hat{\theta}_{N,{j}}^u-\hat{\theta}_{N,{j}}^l)\left[\varphi_j-d_0\right]^+\right\}\\	
            =&{\color{black}\Bigl(\sum_{j=1}^n{\hat{\theta}}^l_{N,j}\Bigr)\mathbb{E}_{P^l}\mathbb{E}_{Q_J}[\varphi(x,\xi)]
+\Bigl(1-\sum_{j=1}^n{\hat{\theta}}^l_{N,j}\Bigr)\inf_{d_0}\left\{d_0+
\frac{\sum_{i=1}^n(\hat{\theta}^u_{N,i}-{\hat{\theta}}^l_{N,i})}{1-\sum_{j=1}^n{\hat{\theta}}^l_{N,j}}
\mathbb{E}_{P^{u-l}}\left[\mathbb{E}_{Q_J}[\varphi(x,\xi)]-d_0\right]^+\right\}}\\
=&{\color{black}\Bigl(\sum_{j=1}^n{\hat{\theta}}^l_{N,j}\Bigr)\mathbb{E}_{P^l}\mathbb{E}_{Q_J}[\varphi(x,\xi)]+\Bigl(1-\sum_{j=1}^n{\hat{\theta}}^l_{N,j}\Bigr)\mathrm{CVaR}^{\beta_N}_{P^{u-l}}\mathbb{E}_{Q_J}[\varphi(x,\xi)]},
	\end{align*}
	where $[a]^+:=\max\{a,0\}$, $j^* := \max\{j\in\{1,\ldots,n\} :\varphi_j<d_0\}$ if $\varphi_1<d_0 $, and $j^ * := 0$ otherwise. The last two equations follow from the definitions of $P^l$ and $P^{u-l}$ and the standard Rockafellar-Uryasev CVaR formula in \cite[Theorem 10]{rockafellar2002conditional}, respectively.
$\hfill\blacksquare$

{\color{black}
\begin{remark}
If $\sum_{j=1}^n \hat\theta^l_{N,j}=0$ (resp.\ $\sum_{j=1}^n(\hat\theta^u_{N,j}-\hat\theta^l_{N,j})=0$), then the corresponding term in the representation vanishes and the distribution
$P^l$ (resp.\ $P^{u-l}$) can be chosen arbitrarily. In particular, if no truncation is active, i.e.,
$\hat\theta_{N,j}-\hat r_N\ge 0$ and $\hat\theta_{N,j}+\hat r_N\le 1$ for all $j$,
we have $\hat\theta^u_{N,j}-\hat\theta^l_{N,j}=2\hat r_N$ for all $j$, and therefore
$P^{u-l}(J=j)=1/n$ (uniform over $\{1,\ldots,n\}$) and $\beta_N=1/2$.
\end{remark}
}

{\color{black}Theorem~\ref{thm-reformulation} reveals that the original lower-level worst-case objective induced by the Bayesian ambiguity set can be expressed as a weighted combination of two kinds of Bayesian risk measures ($\mathbb{E}-\mathbb{E}$ and $\mathrm{CVaR}-\mathbb{E}$), which are widely used in stochastic optimization and control \cite{shapiro2025episodic,wu2018bayesian}. This connects our model to the general relation between DRO and risk-averse optimization. However, in \cite{shapiro2025episodic,wu2018bayesian}, the posterior distribution serves as the outer distribution, and one then performs a prescribed posterior risk functional. In our model, by contrast, the posterior distribution is first used to construct the Bayesian ambiguity set, and the resulting risk representation is then induced endogenously by the geometry of that ambiguity set. Accordingly, the outer measures $P^l$ and $P^{u-l}$, as well as the associated risk level $\beta_N$, are not specified a priori, but arise from the ambiguity-set construction.
}
We next show that
$\nu$ is Lipschitz continuous.

\begin{lemma}\label{lem-con-lip1}
Suppose that \(\mathbb{E}_{Q_j} [\varphi(\cdot, \xi)] \) is Lipschitz continuous with modulus \( L_{\varphi,j} > 0 \) for $j\in\{1,2,\ldots,n\}$. Then $\nu$ is Lipschitz continuous.
\end{lemma}
\textit{Proof.}
Let $(x,\theta,r)$ and $(x',\theta',r')$ be two points in $\mathcal{X}\times\Theta\times[0,1]$.
Denote
\[
\begin{aligned}
\Delta := |\nu(x, \theta, r) - \nu(x', \theta', r')|
\leq \Delta_1 + \Delta_2,
\end{aligned}
\]
where $\Delta_1:=\big|\nu(x,\theta,r)-\nu(x',\theta,r)\big|,
$ and $
\Delta_2:=\big|\nu(x',\theta,r)-\nu(x',\theta',r')\big|.$

\smallskip
\noindent\emph{(i) Bound $\Delta_1$.} Since $\mathbb B(\theta,r)\cap\Theta\subseteq\Theta$ and $\Theta$ is the simplex, we have
\[
\begin{aligned}
   \Delta_1
&=\Big|\max_{\theta^\dagger\in \mathbb B(\theta,r)\cap\Theta}\sum_{j=1}^n \theta_j^\dagger \mathbb E_{Q_j}[\varphi(x,\xi)]
-\max_{\theta^\dagger\in \mathbb B(\theta,r)\cap\Theta}\sum_{j=1}^n \theta_j^\dagger \mathbb E_{Q_j}[\varphi(x',\xi)]\Big|
\\
&\le \max_{\theta^\dagger\in \mathbb B(\theta,r)\cap\Theta}\sum_{j=1}^n \theta_j^\dagger\Big|\mathbb E_{Q_j}[\varphi(x,\xi)]-\mathbb E_{Q_j}[\varphi(x',\xi)]\Big|\le \max_j L_{\varphi,j}\|x-x'\|.
\end{aligned}
\]

\smallskip
\noindent\emph{(ii) Bound $\Delta_2$.}
Define $\Theta':=\mathbb{B}(\theta',r')\cap\Theta$, we first show
\begin{equation}\label{eq:Hbound-inline}
\mathbb H(\mathbb B(\theta,r)\cap\Theta,\Theta')\le \|\theta-\theta'\|_\infty+|r-r'|.
\end{equation}
By definition, it suffices to prove
$\mathbb D(\mathbb B(\theta,r)\cap\Theta,\Theta')\le \|\theta-\theta'\|_\infty+|r-r'|$.
If $\mathbb B(\theta,r)\cap\Theta\subseteq \Theta'$, then $\mathbb D(\mathbb B(\theta,r)\cap\Theta,\Theta')=0$ and the bound holds.
Otherwise, take $\bar\theta\in \mathbb B(\theta,r)\cap\Theta\setminus \Theta'$ and define $\zeta:=\frac{r'}{ \|\bar\theta-\theta'\|_\infty}$. Since $\bar\theta \notin \Theta'$, we have ${r'}<{ \|\bar\theta-\theta'\|_\infty}$, implying that ${\zeta} \in(0,1)$. Define $\tilde{\theta}:={\zeta}\bar\theta+(1-{\zeta})\theta'$. Then,
$$
\mathbf{d}(\tilde{\theta}, \theta')=\|{\zeta} \bar\theta+(1-{\zeta}) \theta'-\theta' \|_\infty\leq {\zeta}\| \bar\theta-\theta' \|_\infty=r',
$$
which means that $\tilde{\theta} \in \Theta'$. Hence
$$
\begin{aligned}
\mathbf{d}\left(\bar\theta, \Theta'\right) & \leq \|\bar\theta-\tilde{\theta}\|_\infty=\|\bar\theta-{\zeta}\bar\theta-(1-{\zeta})\theta'\|_\infty \\
& \leq(1-{\zeta}) \|\bar\theta-\theta'\|_\infty = \|\bar\theta-\theta'\|_\infty -r' \\
& \leq \|\bar\theta-\theta\|_\infty + \|\theta-\theta'\|_\infty -r'\leq r+ \|\theta-\theta'\|_\infty-r'\\
&\leq   \|\theta-\theta'\|_\infty+|r-r'|.
\end{aligned}
$$
By symmetry, the same bound holds for $\mathbb{D}\left(\Theta',\mathbb B(\theta,r)\cap\Theta\right)$, proving \eqref{eq:Hbound-inline}.

Now fix $x'$. Since $\sup_{x\in\mathcal X,j}|\mathbb E_{Q_j}[\varphi(x,\xi)]|\le M$.
Let $\theta^\star\in\argmax_{\theta^\dagger\in\mathbb B(\theta,r)\cap\Theta}\sum_{j=1}^n \theta_j^\dagger\mathbb E_{Q_j}[\varphi(x',\xi)]$.
For any $\varepsilon>0$, by the definition of $\mathbb D(\mathbb B(\theta,r)\cap\Theta,\Theta')$ there exists $\bar\theta'\in\Theta'$ such that
$\|\theta^\star-\bar\theta'\|_\infty\le \mathbb D(\mathbb B(\theta,r)\cap\Theta,\Theta')+\varepsilon$.
Then
\[
\max_{\theta^\dagger\in\mathbb B(\theta,r)\cap\Theta}\sum_{j=1}^n \theta_j^\dagger\mathbb E_{Q_j}[\varphi(x',\xi)]
-\max_{\theta^\dagger\in\Theta'}\sum_{j=1}^n \theta_j^\dagger\mathbb E_{Q_j}[\varphi(x',\xi)]\leq nM\|\theta^\star-\bar\theta'\|_\infty
\le nM\big(\mathbb D(\mathbb B(\theta,r)\cap\Theta,\Theta')+\varepsilon\big).
\]
Letting $\varepsilon\downarrow0$ and using symmetry yields
\[
\Delta_2\le nM\mathbb H(\mathbb B(\theta,r)\cap\Theta,\Theta')
\le nM\big(\|\theta-\theta'\|_\infty+|r-r'|\big).
\]
Combining the bounds on $\Delta_1$ and $\Delta_2$ yields the desired Lipschitz continuity.
$\hfill\blacksquare$

To quantify the discrepancy between the $x$-solution sets $\mathfrak{X}^\lambda(\hat{\theta}_N,\hat{r}_N)$ and $\mathfrak{X}(\theta^c,r^c)$, we first demonstrate that $\mathfrak{X}( \hat{\theta}_N,\hat{r}_N)$ closely approximates ${\mathfrak{X}(\theta^c,r^c)}$. To this end, we present the following lemma.

\begin{lemma}\label{lem-con-svi}
 Suppose that the conditions in Lemma \ref{lem-con-ambiguity} hold and \(\mathbb{E}_{Q_j} [\Phi(\cdot, \xi)] \) is continuous for $j\in\{1,\ldots,n\}$. Then $\mathbb{D}\left(\mathfrak{X}( \hat{\theta}_N,0),{\mathfrak{X}(\theta^c,0)}\right)\to0$ almost surely as $N\to\infty$.
\end{lemma}
 \textit{Proof}.
 Recall that $\mathfrak{X}( \hat{\theta}_N,0)$ and $\mathfrak{X}(\theta^c,0)$ can be regarded as the $x$-solution sets of the SVI problem. Since $\hat{\theta}_N\to\theta^c$ according to Lemma \ref{lem-con-ambiguity}, the result follows directly from Lemma 2.1 in \cite{xu2010sample}.
$\hfill\blacksquare$

\begin{theorem}
Suppose that Assumptions \ref{ass-bvm}-\ref{ass3.1} hold, $\mathbb E_{Q_j}[\Phi(\cdot,\xi)]$ is continuous, and \(\mathbb{E}_{Q_j} [\varphi(\cdot, \xi)] \) is Lipschitz continuous with modulus \( L_{\varphi,j} > 0 \) for $j\in\{1,2,\ldots,n\}$. Then $\mathbb{D}\left(\mathfrak{X}( \hat{\theta}_N,\hat{r}_N),\mathfrak{X}(\theta^c,r^c)\right)\to0$ almost surely as $N\to\infty.$ \label{thm-con-dr}
\end{theorem}
 \textit{Proof}.
Let $\{x_N\}$ be a sequence with $x_N \in \mathfrak{X}( \hat{\theta}_N,\hat{r}_N)$. From the compactness of $\mathcal{X}$,
there exists a subsequence $\{x_{N_k}\}$ of $\{x_N\}$ such that $x_{N_k}\to x^*$ as $k \to \infty$. It suffices to demonstrate that $x^* \in \mathfrak{X}(\theta^c,r^c)$.
For each $x_{N_k} \in \mathfrak{X}( \hat{\theta}_{N_k},\hat{r}_{N_k})$, there exists ${\theta}_{N_k} \in \hat{\Theta}_{N_k}$ satisfying
\begin{equation*}
    \begin{aligned}
        	0 &\in \sum_{j=1}^n{\theta}_{{N_k},j}\mathbb{E}_{Q_j}[\Phi(x_{N_k},\xi)]+\mathcal{N}_\mathcal{X}(x_{N_k}),\\
         {\theta}_{N_k}&\in\underset{\theta\in\hat{\Theta}_{N_k}}{\argmax} \sum_{j=1}^n\theta_j\mathbb{E}_{Q_j}[\varphi\left(x_{N_k}, \xi\right)].
    \end{aligned}
\end{equation*}
Since $\Theta$ is compact, there exists a subsequence (not relabeled) $\{{\theta}_{N_k}\}$ such that  ${\theta}_{N_k}\to{\theta}^*$.
Thus
$$
\begin{aligned}
& \left|\sum_{j=1}^n{\theta}_{N_k,j}\mathbb{E}_{Q_j}[\varphi\left(x_{N_k}, \xi\right)]-\sum_{j=1}^n\theta_{j}^*\mathbb{E}_{Q_j}[\varphi\left(x^*, \xi\right)]\right| \\
\leq & \left|\sum_{j=1}^n{\theta}_{N_k,j}\mathbb{E}_{Q_j}[\varphi\left(x_{N_k}, \xi\right)]-\sum_{j=1}^n{\theta}_{N_k,j}\mathbb{E}_{Q_j}[\varphi\left(x^*, \xi\right)]\right| +\left|\sum_{j=1}^n{\theta}_{N_k,j}\mathbb{E}_{Q_j}[\varphi\left(x^*, \xi\right)]-\sum_{j=1}^n{\theta}_{j}^*\mathbb{E}_{Q_j}[\varphi\left(x^*, \xi\right)]\right| \\
\leq & \sum_{j=1}^n{\theta}_{N_k,j}L_{\varphi,j}\left\|x_{N_k}-x^*\right\|+\sum_{j=1}^n|{\theta}^*_{j}-{\theta}_{N_k,j}|\left|\mathbb{E}_{Q_j}[\varphi\left(x^*, \xi\right)] \right| 
\to  0,
\end{aligned}
$$
as $k \to \infty$. Based on the convergence in Theorem \ref{thm-new} and the Lipschitz property of $\nu$ in Lemma \ref{lem-con-lip1}, it follows that $\nu(x_{N_k},  \hat{\theta}_{N_k},\hat{r}_{N_k})\to\nu(x^*, \theta^c,r^c)$ almost surely as $k \to \infty$. Since $\sum_{j=1}^n{\theta}_{N_k,j}\mathbb{E}_{Q_j}[\varphi\left(x_{N_k}, \xi\right)]= \nu(x_{N_k},  \hat{\theta}_{N_k},\hat{r}_{N_k})$, we have $\sum_{j=1}^n{\theta}_{j}^*\mathbb{E}_{Q_j}[\varphi\left(x^*, \xi\right)]=\nu(x^*,  \theta^c,r^c)$, which leads to
$
{\theta}^*\in \vartheta(x^*,  \theta^c,r^c).
$

As $k \to \infty$, we have $\sum_{j=1}^n {\theta}_{N_k,j} \mathbb{E}_{Q_j}[\Phi\left(x_{N_k}, \xi\right)] \to \sum_{j=1}^n{\theta}_{j}^* \mathbb{E}_{Q_j}[\Phi\left(x^*, \xi\right)]$. Since $\mathcal{N}_{\mathcal{X}}(\cdot)$ is upper semicontinuous for the closed and convex $\mathcal{X}$ and
$
0 \in \sum_{j=1}^n {\theta}_{N_k,j} \mathbb{E}_{Q_j}[\Phi\left(x_{N_k}, \xi\right)]+\mathcal{N}_{\mathcal{X}}\left(x_{N_k}\right),
$
we obtain
$
0 \in\sum_{j=1}^n {\theta}_{j}^* \mathbb{E}_{Q_j}[\Phi\left(x^*, \xi\right)]+\mathcal{N}_{\mathcal{X}}\left(x^*\right),
$
which completes the proof.
$\hfill\blacksquare$

\begin{remark}
When we choose $r^c = 0$, Theorem~\ref{thm-con-dr} yields the convergence result
$\mathbb{D}\left(\mathfrak{X}(\hat{\theta}_N, \hat{r}_N), \mathfrak{X}(\theta^c, 0)\right) \to 0,$
i.e., in the risk-neutral case, the Bayesian DRVI solutions converge to a solution of SVI \eqref{SVI} with the true distribution $P_{\theta^c}$. This means that, in practice, our ambiguity-aware decisions become increasingly accurate and less conservative as parametric uncertainty diminishes with the accumulation of information or data.
\end{remark}

In the following, for a fixed sample size $N$ with the ambiguity set $\hat{\Theta}_N$, we analyze the convergence of the $x$-solution set $\mathfrak{X}^\lambda\left(\hat{\theta}_N, \hat{r}_N\right)$ of \eqref{bdrvi-regular2} as $\lambda$ tends to 0. We begin by establishing the uniform convergence of the regularized optimal value function to the corresponding optimal value function of problem \eqref{bdrvi2}.

\begin{lemma}\label{lem-con-regular}
For any fixed $\hat{\theta}_N,\hat{r}_N$, the optimal value $\nu^\lambda(\cdot,  \hat{\theta}_N,\hat{r}_N)$
converges to $\nu(\cdot,  \hat{\theta}_N,\hat{r}_N)$ uniformly as $\lambda\downarrow0$ over  $x\in\mathcal{X}$.
Moreover, for any fixed $x \in \mathcal{X}$, the unique maximizer of the regularized lower-level problem converges to the unique least-norm solution of the unregularized problem, i.e.,
\[
\vartheta^\lambda(x,\hat{\theta}_N,\hat{r}_N)\to \theta_N^{\min}
:=\argmin_{\theta\in\vartheta(x,\hat{\theta}_N,\hat{r}_N)}\|\theta\|.
\]
\end{lemma}
 \textit{Proof}.
We begin by proving the uniform convergence of $\nu^\lambda(\cdot,  \hat{\theta}_N,\hat{r}_N)$ to $\nu(\cdot,  \hat{\theta}_N,\hat{r}_N)$. {\color{black}For any fixed $x\in\mathcal X$, the functions
$
\theta\mapsto \sum_{j=1}^n\theta_j\mathbb E_{Q_j}[\varphi(x,\xi)]
\quad\text{and}\quad
\theta\mapsto \sum_{j=1}^n\theta_j\mathbb E_{Q_j}[\varphi(x,\xi)]-\lambda\|\theta\|^2
$
are continuous on the compact set $\hat{\Theta}_N$.
Hence the optimal values $\nu(x,\hat\theta_N,\hat r_N)$ and
$\nu^\lambda(x,\hat\theta_N,\hat r_N)$ are well defined.}
For any given $\lambda>0$, it is evident that
$
    \nu^\lambda(x,  \hat{\theta}_N,\hat{r}_N) \leq \nu(x,  \hat{\theta}_N,\hat{r}_N).
$
Let
$\theta^\lambda :=\vartheta^\lambda(x,  \hat{\theta}_N,\hat{r}_N)$, since $\|\theta\|^2\leq1$ for all $\theta\in\hat{\Theta}_N$,
we have that
\begin{equation}
   \nu^\lambda(x,  \hat{\theta}_N,\hat{r}_N) =\sum_{j=1}^n\theta_j^\lambda\mathbb{E}_{Q_j}[\varphi(x,\xi)]- \lambda \|\theta^\lambda\|^2=\max_{\theta\in\hat{\Theta}_N} \sum_{j=1}^n\theta_j\mathbb{E}_{Q_j}[\varphi(x,\xi)]- \lambda \|\theta\|^2\geq\nu(x, \hat{\theta}_N, \hat{r}_N)- \lambda,
\end{equation}
which means $|\nu^\lambda(x,  \hat{\theta}_N,\hat{r}_N)-\nu(x,  \hat{\theta}_N,\hat{r}_N)|\leq\lambda$ for all $x\in\mathcal{X}$.
Thus, we can obtain that
$\nu^{\lambda}(\cdot,  \hat{\theta}_N,\hat{r}_N)$ converges to $\nu(\cdot,  \hat{\theta}_N,\hat{r}_N)$ uniformly as $\lambda\downarrow0$.

We now turn to the second claim. For any fixed $x\in\mathcal X$, the set $\vartheta(x,  \hat{\theta}_N,\hat{r}_N)$ is nonempty compact and convex. Since $\|\cdot\|^2$ is strictly convex, there exists a unique $\theta_N^{\min}$ such that $\theta_N^{\min}=\argmin\{\|\theta\|:\theta\in \vartheta(x,  \hat{\theta}_N,\hat{r}_N)\}$.
Consider any sequence $\{\lambda_k\}$ with $\lambda_k\downarrow0$. Since $\theta^{\lambda_k}\in\hat{\Theta}_N$ and $\hat{\Theta}_N$ is compact, there exists a subsequence (not relabeled) $\{\theta^{\lambda_k}\}$ such that $\theta^{\lambda_k}\to\bar{\theta}\in\hat{\Theta}_N$.
By the optimality of $\theta^{\lambda_k}$, for any $\hat{\theta}\in\vartheta(x,\hat{\theta}_N,\hat{r}_N)$,
\[
\sum_{j=1}^n\theta^{\lambda_k}_j\mathbb{E}_{Q_j}[\varphi(x,\xi)]-\lambda_k\|\theta^{\lambda_k}\|^2
\ \ge\
\sum_{j=1}^n\hat{\theta}_j\mathbb{E}_{Q_j}[\varphi(x,\xi)]-\lambda_k\|\hat{\theta}\|^2.
\]

Taking the limit as $k\to\infty$, we obtain
$\sum_{j=1}^n\bar{\theta}_j\mathbb{E}_{Q_j}[\varphi(x,\xi)]\ \ge\ \sum_{j=1}^n\hat{\theta}_j\mathbb{E}_{Q_j}[\varphi(x,\xi)]$ for all $\hat{\theta}\in\vartheta(x,\hat{\theta}_N,\hat{r}_N)$,
which implies $\bar{\theta}\in\vartheta(x,\hat{\theta}_N,\hat{r}_N)$.
Taking $\hat{\theta}=\theta^{\min}_N$ leads to
\[
\lambda_k\bigl(\|\theta^{\lambda_k}\|^2-\|\theta^{\min}_N\|^2\bigr)
\leq
\sum_{j=1}^n\theta^{\lambda_k}_j\mathbb{E}_{Q_j}[\varphi(x,\xi)]-\sum_{j=1}^n\theta^{\min}_{N,j}\mathbb{E}_{Q_j}[\varphi(x,\xi)]
\leq 0,
\]
so $\|\theta^{\lambda_k}\|\le \|\theta^{\min}_N\|$ for all $k$, and therefore
$
\|\bar\theta\|\leq\liminf_{k\to\infty}\|\theta^{\lambda_k}\|\leq\|\theta^{\min}_N\|.
$
By the uniqueness of the least-norm element, it follows that $\bar\theta=\theta^{\min}_N$. Since every sequence $\{\theta^{\lambda_k}\}$ with $\lambda_k\downarrow0$ admits a subsequence converging to $\theta^{\min}_N$, we conclude that $\theta^\lambda \to \theta^{\min}_N$ as $\lambda \downarrow 0$.
$\hfill\blacksquare$

\begin{theorem}
Suppose that \(\mathbb{E}_{Q_j} [\Phi(\cdot, \xi)] \) and \(\mathbb{E}_{Q_j} [\varphi(\cdot, \xi)] \) are Lipschitz continuous with modulus \( L_{\Phi,j} > 0 \) and \( L_{\varphi,j} > 0 \) for $j\in\{1,2,\ldots,n\}$, respectively. Then $\mathbb{D}\left(\mathfrak{X}^\lambda( \hat{\theta}_N,\hat{r}_N),\mathfrak{X}( \hat{\theta}_N,\hat{r}_N)\right)\to0$ as $\lambda\downarrow0.$ \label{thm-con-dr-regular}
\end{theorem}
 \textit{Proof}.
Consider a non-increasing sequence $\{\lambda_k\}$ that converges to 0 as $k \to \infty$ with $x^k_N \in\mathfrak{X}^{\lambda_k}( \hat{\theta}_N,\hat{r}_N)$. By compactness of $\mathcal X$, there exists a subsequence (not relabeled) $\{x_N^k\}$ with $x^k_N \to x_N^*$ as $k \to \infty$. Our goal is to verify that $x_N^* \in \mathfrak{X}(\hat{\theta}_N,\hat{r}_N)$. By the definition of $\mathfrak{X}^{\lambda_k}( \hat{\theta}_N,\hat{r}_N)$, for each $x^k_N$, there exists ${\theta}^k_{N} \in \hat{\Theta}_N$ such that
\begin{equation*}
    \begin{aligned}
        	0 &\in \sum_{j=1}^n{\theta}^k_{N,j}\mathbb{E}_{Q_j}[\Phi\left(x^k_N, \xi\right)]+\mathcal{N}_\mathcal{X}(x^k_N),\\
         {\theta}^k_{N}&\in\underset{\theta\in\hat{\Theta}_N}{\argmax} \sum_{j=1}^n\theta_j\mathbb{E}_{Q_j}[\varphi\left(x^k_N, \xi\right)]-\lambda_k\|\theta\|^2.
    \end{aligned}
\end{equation*}
Since $\theta_N^k\in\hat \Theta_N$ and $\hat \Theta_N$ is compact, there exists a subsequence (not relabeled) $\{{\theta}^k_{N}\}$ with ${\theta}^k_{N}\to{\theta}_N^*$. Then we have
\[
\begin{aligned}
&\left| \sum_{j=1}^n {\theta}_{N,j}^k \mathbb{E}_{Q_j}[\varphi\left(x^k_N, \xi\right)] - \lambda_k \|{\theta}_N^k\|^2 - \sum_{j=1}^n {\theta}_{N,j}^* \mathbb{E}_{Q_j}[\varphi\left(x^*_N, \xi\right)] \right| \\
&\leq \sum_{j=1}^n |{\theta}_{N,j}^k| \cdot |\mathbb{E}_{Q_j}[\varphi\left(x^k_N, \xi\right)] - \mathbb{E}_{Q_j}[\varphi\left(x^*_N, \xi\right)]|
+ \sum_{j=1}^n |{\theta}_{N,j}^k - {\theta}_{N,j}^*| \cdot |\mathbb{E}_{Q_j}[\varphi\left(x^*_N, \xi\right)]|
+ \lambda_k \|{\theta}_N^k\|^2 \to 0.
\end{aligned}
\]
Meanwhile, by Lemmas \ref{lem-con-lip1} and \ref{lem-con-regular}, we have
$\nu^{\lambda_k}(x_N^k, \hat{\theta}_N, \hat{r}_N) \to \nu(x^*_N, \hat{\theta}_N, \hat{r}_N)$.
With this and the fact that $\nu^{\lambda_k}(x_N^k, \hat{\theta}_N, \hat{r}_N) = \sum_{j=1}^n {\theta}_{N,j}^k \mathbb{E}_{Q_j}[\varphi\left(x^k_N, \xi\right)] - \lambda_k \|{\theta}_N^k\|^2$, we can conclude that
$\sum_{j=1}^n {\theta}_{N,j}^* \mathbb{E}_{Q_j}[\varphi\left(x^*_N, \xi\right)] = \nu(x_N^*, \hat{\theta}_N, \hat{r}_N)$, which leads to $ {\theta}_N^* \in \vartheta(x_N^*, \hat{\theta}_N, \hat{r}_N).$

Following a procedure similar to the proof of Theorem \ref{thm-con-dr},
along with
$0 \in \sum_{j=1}^n {\theta}^k_{N,j}\mathbb{E}_{Q_j}[ \Phi\left(x^k_N, \xi\right)]+\mathcal{N}_{\mathcal{X}}\left(x_N^k\right)$,
we can obtain that
$
0 \in\sum_{j=1}^n {\theta}_{N,j}^* \mathbb{E}_{Q_j}[\Phi\left(x^*_N, \xi\right)]+\mathcal{N}_{\mathcal{X}}\left(x^*_N\right),
$
which completes the proof.
$\hfill\blacksquare$

\begin{remark}\label{cor:theta-mn-conv}
Under the conditions in Theorem~\ref{thm-con-dr-regular}, there exists a sequence $\{\lambda_k\}$ with $\lambda_k\downarrow 0$ such that
$x_N^k\in\mathfrak X^{\lambda_k}(\hat{\theta}_N,\hat r_N)$ with $x_N^k\to x_N^\ast$.
For $\theta_N^k\in\vartheta^{\lambda_k}(x_N^k,\hat{\theta}_N,\hat r_N)$,
there exists a subsequence such that $\theta_N^k\to\theta_N^\ast$
and $(x_N^\ast,\theta_N^\ast)$ is a solution of \eqref{bdrvi2}. Moreover, from Lemma \ref{lem-con-regular}, $\theta_N^\ast$ is the unique least-norm optimal solution of $\vartheta(x_N^\ast,\hat{\theta}_N,\hat r_N)$, i.e.,
$\theta_N^\ast\equiv\theta_N^{\min}=\argmin\Bigl\{\|\theta\|:\ \theta\in \vartheta(x_N^\ast,\hat{\theta}_N,\hat r_N)\Bigr\}.$
\end{remark}

We are now ready to present the main convergence theorem for this section.
\begin{theorem}\label{thm-con}
      Suppose that Assumptions \ref{ass-bvm}-\ref{ass3.1} and the conditions in Theorem \ref{thm-con-dr-regular} hold. Then, almost surely, $$\lim_{N\to\infty}\lim_{\lambda\downarrow0}\mathbb{D}\left(\mathfrak{X}^{\lambda}( \hat{\theta}_N,\hat{r}_N),{\mathfrak{X}(\theta^c,r^c)}\right)=0.$$
\end{theorem}
 \textit{Proof}.
Based on the triangle inequality, we have
$$
\begin{aligned}
\mathbb{D}\left(\mathfrak{X}^\lambda\left(\hat{\theta}_N, \hat{r}_N\right), \mathfrak{X}\left(\theta^c,r^c\right)\right)
&\leq \mathbb{D}\left(\mathfrak{X}^\lambda\left(\hat{\theta}_N, \hat{r}_N\right), \mathfrak{X}\left(\hat{\theta}_N, \hat{r}_N\right)\right)+\mathbb{D}\left(\mathfrak{X}\left(\hat{\theta}_N, \hat{r}_N\right), \mathfrak{X}\left(\theta^c,r^c\right)\right).
\end{aligned}
$$
The first term vanishes as $\lambda\downarrow0$ by Theorem~\ref{thm-con-dr-regular}, and the second term vanishes almost surely as $N\to\infty$ by Theorem \ref{thm-con-dr}. Hence the left-hand side converges to zero almost surely as claimed.
$\hfill\blacksquare$

{\color{black}
Theorem~\ref{thm-con} provides a sequential asymptotic consistency result, obtained by first letting $\lambda\downarrow0$ for fixed $N$ and then letting $N\to\infty$. While this ensures convergence to the Bayesian DRVI solution set, it does not yet determine how $\lambda$ should be chosen as a function of $N$. In the next section, we address this issue through quantitative stability bounds, which make it possible to discuss the trade-off between sampling error and regularization bias.
}

\section{Stability of solution sets}
While Section 3 proves a sequential asymptotic  consistency result for solutions to (\ref{bdrvi-regular2}) as $N\to\infty$ and $\lambda\downarrow0$, practical applications rely on finite, often noisy datasets and a fixed positive $\lambda$.
This setting introduces two distinct sources of deviation from the idealized solution:
(i) sampling noise, which perturbs the estimated ambiguity parameters $(\hat{\theta}_N,\hat r_N)$, and (ii) the regularization bias due to $\lambda>0$.
These facts highlight the importance of stability analysis in ensuring robustness.
In particular, small perturbations in $\mathcal{S}^N$ may generate a contaminated dataset $\tilde{\mathcal{S}}^N=\{\tilde{\xi}^1,\dots,\tilde{\xi}^N\}$ with shifted ambiguity parameters $(\tilde{\theta}_N,\tilde{r}_N)$, leading to the perturbed Bayesian DRVI problem:
\begin{equation}
    \begin{aligned}
        	0 &\in \sum_{j=1}^n\theta_j\mathbb{E}_{Q_j}[\Phi(x,\xi)]+\mathcal{N}_\mathcal{X}(x),\\
         \theta&\in\underset{{\theta}^\prime\in\tilde{\Theta}_N}{\argmax} \sum_{j=1}^n{\theta}^\prime_j\mathbb{E}_{Q_j}[\varphi(x,\xi)]-\lambda\|{\theta}^\prime\|^2,
    \end{aligned}\label{p-bdrvi2}\tag{$\tilde{\text{B}}_N^\lambda$}
\end{equation}
where the perturbed ambiguity set $\tilde{\Theta}_N:=\mathbb{B}(\tilde{\theta}_N,\tilde{r}_N)\cap\Theta$ is constructed from the contaminated dataset ${\tilde{\mathcal{S}}^{N}}$.

This section investigates how the $x$-solution set responds to perturbations of ambiguity sets. In addition, we provide finite-sample guarantees that ensure statistical performance with high probability. Further, we analyze whether the perturbed solution set $\mathfrak{X}^\lambda(\tilde{\theta}_N,\tilde{r}_N)$ closely approximates the solution set $\mathfrak{X}(\theta^c,r^c)$ for sufficiently large $N$. If such stability holds, the perturbed solutions can be justified as reliable approximations to \eqref{bdrvi-true} even under data contamination.
We begin by establishing the following lemma, which provides a bound on the Hausdorff distance between ambiguity sets corresponding to different parameters.
\begin{lemma}[Quantitative stability of Bayesian parametric ambiguity sets]\label{lem-q-BAS}
For any $\hat{\theta}_N$, $\tilde{\theta}_N \in \Theta$ and $\hat{r}_N, \tilde{r}_N \in[0,1]$, we have
$$
\mathbb{H}\left(\hat{\Theta}_N,\tilde{\Theta}_N\right) \leq \|\hat{\theta}_N-\tilde{\theta}_N\|_\infty+\left|\tilde{r}_N-\hat{r}_N\right|.
$$
\end{lemma}
This is a direct consequence by the proof of Lemma~\ref{lem-con-lip1} from \eqref{eq:Hbound-inline}, we omit the proof.

\begin{lemma}\label{lem-lip-optimal}
Suppose that the conditions in Theorem \ref{thm-con} hold. Then, for any $x\in\mathcal{X}$ and $\lambda\geq0$,
$$|\nu^\lambda(x, \hat{\theta}_N,\hat{r}_N)-\nu^\lambda(x,\tilde{\theta}_N,\tilde{r}_N)| \leq  (nM+2\lambda)(\|\hat{\theta}_N-\tilde{\theta}_N\|_\infty+\left|\tilde{r}_N-\hat{r}_N\right|).$$
\end{lemma}
 \textit{Proof}.
 Note that
$$
\begin{aligned}
&\nu^\lambda(x, \hat{\theta}_N,\hat{r}_N)-\nu^\lambda(x,\tilde{\theta}_N,\tilde{r}_N)\\ & =\max _{\theta \in \hat{\Theta}_N} \left[\sum_{j=1}^n\theta_j\mathbb{E}_{Q_j}[\varphi(x,\xi)]-\lambda\|\theta\|^2\right] -\max _{\theta' \in \tilde{\Theta}_N} \left[\sum_{j=1}^n\theta_j'\mathbb{E}_{Q_j}[\varphi(x,\xi)]-\lambda\|\theta'\|^2\right]  \\
& \leq\max _{\theta \in \hat{\Theta}_N}\min _{\theta' \in \tilde{\Theta}_N} \left(\sum_{j=1}^n\theta_j\mathbb{E}_{Q_j}[\varphi(x,\xi)]-\lambda\|\theta\|^2-\sum_{j=1}^n\theta_j'\mathbb{E}_{Q_j}[\varphi(x,\xi)]+\lambda\|\theta'\|^2\right) \\
& \leq nM \max_{\theta\in \hat{\Theta}_N}\min _{\theta' \in \tilde{\Theta}_N}\|\theta-\theta'\|_\infty +\lambda\max_{\theta\in \hat{\Theta}_N}\min _{\theta' \in \tilde{\Theta}_N}\left|\|\theta\|^2-\|\theta'\|^2\right| \\
& \leq nM \max_{\theta\in \hat{\Theta}_N}\min _{\theta' \in \tilde{\Theta}_N}\|\theta-\theta'\|_\infty +2\lambda\max_{\theta\in \hat{\Theta}_N}\min _{\theta' \in \tilde{\Theta}_N}\|\theta-\theta'\|_\infty  =(nM+2\lambda)\mathbb{D}\left(\hat{\Theta}_N,\tilde{\Theta}_N\right).
\end{aligned}
$$
Similarly, we can show that $\nu^\lambda(x, \tilde{\theta}_N,\tilde{r}_N)-\nu^\lambda(x,\hat{\theta}_N,\hat{r}_N) \leq (nM+2\lambda)\mathbb{D}\left(\tilde{\Theta}_N,\hat{\Theta}_N\right)$. Thus, combining with Lemma \ref{lem-q-BAS}, we can draw the conclusion.
$\hfill\blacksquare$

The next result quantifies the stability of the lower-level maximizer set
$\vartheta^\lambda(x,\hat{\theta}_N,\hat{r}_N)$ under perturbations of the ambiguity set.
Although the lower-level maximization problem admits a unique solution when $\lambda>0$,
we retain the set-valued notation $\vartheta^\lambda$ in order to cover the case $\lambda=0$ as well.
{\color{black}To derive quantitative stability estimates, it is not enough to know that the lower-level optimal value changes continuously under perturbations.
One also needs a way to relate small value perturbations to small changes in the corresponding optimizer set.}
Inspired by the literature in stochastic optimization \cite{romisch2003stability}, we consider the growth function $\psi_x^\lambda$ defined on $\mathbb{R}_{+}$ by
$$
\psi_x^\lambda(\tau):=\min_{ \theta \in \hat{\Theta}_N} \left\{  \nu^\lambda(x,\hat{\theta}_N,\hat{r}_N)-\left(\sum_{j=1}^n\theta_j\mathbb{E}_{Q_j}[\varphi(x,\xi)]-\lambda\|\theta\|^2 \right): \mathbf{d}\left(\theta,\vartheta^\lambda(x,\hat{\theta}_N,\hat{r}_N)\right) \geq \tau\right\},
$$
together with the associated function
$\Psi_x^\lambda(\eta):=\eta+(\psi_x^\lambda)^{-1}(2 \eta)$,
where $(\psi_x^\lambda)^{-1}(t):=\sup \left\{\tau \in \mathbb{R}_{+}: \psi_x^\lambda(\tau) \leq t\right\}$ and $\eta \in \mathbb{R}_{+}$.
{\color{black}
Intuitively, $\psi_x^\lambda(\tau)$ measures the smallest loss in the lower-level objective incurred by choosing a parameter $\theta$ that stays at least $\tau$ away from the maximizer set $\vartheta^\lambda(x,\hat{\theta}_N,\hat{r}_N)$.
Accordingly, the inverse $(\psi_x^\lambda)^{-1}$ converts lower-level value perturbations into distance bounds for the corresponding worst-case parameter, which is the key mechanism used in the following stability analysis.
}

\begin{theorem}\label{thm-uppercon-1}
Suppose that the conditions in Theorem \ref{thm-con} hold. Then, for any fixed $x\in\mathcal{X}$,
there exists a constant $\hat{L} \geq 1$ such that
$$
\begin{aligned}
    \vartheta^\lambda(x,\tilde{\theta}_N,\tilde{r}_N) &\subseteq \vartheta^\lambda(x,\hat{\theta}_N,\hat{r}_N)+\Psi_x^\lambda\left(\hat{L} \mathbb{D}(\tilde{\Theta}_N,\hat{\Theta}_N)\right) \mathbb{B}_\infty\\&\subseteq \vartheta^\lambda(x,\hat{\theta}_N,\hat{r}_N)+\Psi_x^\lambda\left(\hat{L} (\|\tilde{\theta}_N-\hat{\theta}_N\|_\infty+\left|\tilde{r}_N-\hat{r}_N\right|)\right) \mathbb{B}_\infty
\end{aligned}
$$
holds for any $\tilde{\theta}_N\in\Theta$ and $\tilde{r}_N\in[0,1]$.
Here $\mathbb B_\infty$ denotes the unit ball under the $\ell_\infty$ norm.
\end{theorem}
 \textit{Proof}.   The argument closely follows the structure of \cite[Theorem 9]{romisch2003stability}.
Let ${\theta}^*\in\vartheta^\lambda(x,\tilde{\theta}_N,\tilde{r}_N)$ and $\bar{\delta}:=\mathbb{D}(\tilde{\Theta}_N,\hat{\Theta}_N)$,
there exists $\bar{\theta}\in\hat{\Theta}_N$ such that $\|{\theta}^*-\bar{\theta}\|_\infty \leq \bar{\delta}$.
From Lemma \ref{lem-lip-optimal} and the definition of $\psi_x^\lambda$, we have
$$
\begin{aligned}
2(nM+2\lambda)\bar{\delta} & \geq (nM+2\lambda)\bar{\delta} +\nu^\lambda(x,\hat{\theta}_N,\hat{r}_N)-\nu^\lambda(x,\tilde{\theta}_N,\tilde{r}_N) \\
& = (nM+2\lambda)\bar{\delta}+\nu^\lambda(x,\hat{\theta}_N,\hat{r}_N)-\left(\sum_{j=1}^n{\theta}^*_j\mathbb{E}_{Q_j}[\varphi(x,\xi)]-\lambda\|{\theta}^*\|^2\right) \\
& \geq\nu^\lambda(x,\hat{\theta}_N,\hat{r}_N) -\left(\sum_{j=1}^n\bar{\theta}_j\mathbb{E}_{Q_j}[\varphi(x,\xi)] -\lambda\|\bar{\theta}\|^2\right) \geq \psi_x^\lambda\left(\mathbf{d}\left(\bar{\theta},\vartheta^\lambda(x,\hat{\theta}_N,\hat{r}_N)\right)\right) \\
& \geq \inf _{\theta \in {\theta}^*+\bar{\delta} \mathbb{B}_\infty} \psi_x^\lambda\left(\mathbf{d}\left(\theta, \vartheta^\lambda(x,\hat{\theta}_N,\hat{r}_N)\right)\right)=\psi_x^\lambda\left(\mathbf{d}\left({\theta}^*, \vartheta^\lambda(x,\hat{\theta}_N,\hat{r}_N)+ \bar{\delta} \mathbb{B}_\infty\right)\right),
\end{aligned}
$$
where the last inequality follows from the monotonicity of $\psi_x^\lambda$. This implies
$$
\begin{aligned}
\mathbf{d}\left({\theta}^*, \vartheta^\lambda(x,\hat{\theta}_N,\hat{r}_N)\right) & \leq  \bar{\delta}+\mathbf{d}\left({\theta}^*, \vartheta^\lambda(x,\hat{\theta}_N,\hat{r}_N)+ \bar{\delta} \mathbb{B}_\infty\right)\\
& \leq\bar{\delta}+(\psi_x^\lambda)^{-1}\left(2(nM+2\lambda)\bar{\delta}\right) \leq \hat{L} \bar{\delta}+(\psi_x^\lambda)^{-1}(2 \hat{L} \bar{\delta})=\Psi_x^\lambda(\hat{L} \bar{\delta})
\end{aligned}
$$
where $\hat{L}:=\max \left\{1, nM+2\lambda\right\}$.
Thus we have $\vartheta^\lambda(x,\tilde{\theta}_N,\tilde{r}_N) \subseteq \vartheta^\lambda(x,\hat{\theta}_N,\hat{r}_N)+\Psi_x^\lambda\left(\hat{L} \mathbb{D}(\tilde{\Theta}_N,\hat{\Theta}_N)\right) \mathbb{B}_\infty$.
Since $\Psi_x^\lambda$ is monotonically increasing, we further conclude that
$$
\begin{aligned}
\vartheta^\lambda(x,\hat{\theta}_N,\hat{r}_N)+\Psi_x^\lambda\left(\hat{L} \mathbb{D}( \tilde{\Theta}_N,\hat{\Theta}_N)\right) \mathbb{B}_\infty
\subseteq \vartheta^\lambda(x,\hat{\theta}_N,\hat{r}_N)+\Psi_x^\lambda\left(\hat{L} (\|\tilde{\theta}_N-\hat{\theta}_N\|_\infty+\left|\tilde{r}_N-\hat{r}_N\right|)\right) \mathbb{B}_\infty.
\end{aligned}
$$
Thus, the theorem follows.
$\hfill\blacksquare$

\begin{remark}\label{remark-2}
Suppose that there exist $p_\lambda\geq1$ and $\gamma_\lambda>0$ such that for each sufficiently small $\tau \in \mathbb{R}_{+}$, we have $\psi_x^\lambda(\tau)\geq\gamma_\lambda\tau^{p_\lambda}.$
Then for each $\theta$ close to $\vartheta^\lambda(x,\hat{\theta}_N,\hat{r}_N)$, the following holds:
$$
\nu^\lambda(x,\hat{\theta}_N,\hat{r}_N) \geq\sum_{j=1}^n\theta_j\mathbb{E}_{Q_j}[\varphi(x,\xi)] -\lambda\|\theta\|^2+\gamma_\lambda\mathbf{d}\left({\theta}, \vartheta^\lambda(x,\hat{\theta}_N,\hat{r}_N)\right)^{p_\lambda}.$$
Under these conditions, the function $\Psi_x^\lambda(\eta)$ can be bounded as $\Psi_x^\lambda(\eta)\leq \eta+\left(\frac{2 \eta}{\gamma_\lambda}\right)^{1/p_\lambda} \leq C \eta^{1/p_\lambda}$ for some constant $C>0$ and sufficiently small $\eta \in \mathbb{R}_{+}$.
\end{remark}

\begin{remark} \label{remark-3}
When \(\lambda > 0\), the regularized objective function \(\sum_{j=1}^n \theta_j \mathbb{E}_{Q_j}[\varphi(x, \xi)] - \lambda \|\theta\|^2\) is strictly concave. Near the unique maximizer \(\theta^* \in \hat{\Theta}_N\),  a second-order expansion gives:
\[
\nu^\lambda(x, \hat{\theta}_N, \hat{r}_N) - \left( \sum_{j=1}^n \theta_j \mathbb{E}_{Q_j}[\varphi(x, \xi)] - \lambda \|\theta\|^2 \right) \geq \lambda \|\theta - \theta^*\|^2 \geq c \lambda \|\theta - \theta^*\|_\infty^2
\]
for some $c > 0$ by norm equivalence.
Thus, the growth function satisfies \( \psi_x^\lambda(\tau) \geq c\lambda \tau^2 \).
For \(\lambda = 0\), the objective function \(\sum_{j=1}^n \theta_j \mathbb{E}_{Q_j}[\varphi(x, \xi)]\) becomes linear in \(\theta\). Define
$\delta(x):=\max_j \mathbb{E}_{Q_j}[\varphi(x, \xi)]-\max_{j\notin I(x)}\mathbb{E}_{Q_j}[\varphi(x, \xi)]$ and $I(x):=\argmax_j \mathbb{E}_{Q_j}[\varphi(x, \xi)]$.
Assuming that there exists
$\underline{\delta}:=\inf_{x\in\mathcal{X}}\delta(x)>0, $
then the growth function satisfies \( \psi_x(\tau) \geq \frac{1}{2}\underline{\delta} \tau \).
\end{remark}

{\color{black}
When $\lambda=0$, the lower-level problem is linear over $\hat{\Theta}_N$ and may therefore admit multiple maximizers, so first-order growth generally requires an additional nondegeneracy condition. By contrast, when $\lambda>0$, the term $-\lambda\|\theta\|^2$ makes the lower-level objective strongly concave in $\theta$, yields a unique maximizer, and automatically guarantees second-order growth. 
}
As a direct consequence of Remarks \ref{remark-2} and \ref{remark-3}, Proposition \ref{prop-H} establishes that the mapping \( \vartheta^\lambda(x,\cdot,\cdot) \) is Hölder continuous at \( (\hat{\theta}_N,\hat{r}_N) \) with exponent \( \frac{1}{p_\lambda} \) in the regularized-to-regularized case, while the regularized-to-unregularized deviation is controlled with exponent \(\frac{1}{p_0} \).
\begin{proposition}\label{prop-H}
Suppose that the conditions in Theorem \ref{thm-con} hold.  Assume that for any fixed $\lambda\geq0$, there exist $p_\lambda\geq1$ and $\gamma_\lambda>0$ such that $\psi_x^\lambda(\tau)\geq\gamma_\lambda\tau^{p_\lambda}$  for all $x\in\mathcal{X}$.
Then
for any $\hat{\theta}_N$, $\tilde{\theta}_N \in \Theta$ and $\hat{r}_N, \tilde{r}_N \in [0,1]$,
we have
$$
\sup_{x\in\mathcal{X}}\mathbb{H}\left(\vartheta^\lambda(x,\tilde{\theta}_N, \tilde{r}_N),\vartheta^\lambda(x,\hat{\theta}_N, \hat{r}_N)\right) \leq \left(\frac{1}{\gamma_\lambda}(nM+2\lambda)(\|\hat{\theta}_N-\tilde{\theta}_N\|_\infty+\left|\tilde{r}_N-\hat{r}_N\right|)\right)^{1/p_\lambda},$$
$$
\sup_{x\in\mathcal{X}}\mathbb{D}\left(\vartheta^\lambda(x,\tilde{\theta}_N, \tilde{r}_N),\vartheta(x,\hat{\theta}_N, \hat{r}_N)\right) \leq \left(\frac{1}{\gamma_0}(nM+2\lambda)(\|\hat{\theta}_N-\tilde{\theta}_N\|_\infty+\left|\tilde{r}_N-\hat{r}_N\right|)+\frac{\lambda}{\gamma_0}\right)^{1/p_0}.$$
\end{proposition}
 \textit{Proof}.
For the first part, let $\tilde{\theta}^*\in\vartheta^\lambda(x,\tilde{\theta}_N, \tilde{r}_N)$ and $\theta^*\in\vartheta^\lambda(x,\hat{\theta}_N, \hat{r}_N)$ be such that $\|\theta^*-\tilde{\theta}^*\|_\infty=\mathbf{d}\left(\tilde{\theta}^*, \vartheta^\lambda(x,\hat{\theta}_N, \hat{r}_N)\right)$. This leads to the following estimation:
\begin{equation*}
    \begin{aligned}
        \psi_x^\lambda(\|\theta^*-\tilde{\theta}^*\|_\infty)&=\psi_x^\lambda\left(\mathbf{d}\left(\tilde{\theta}^*, \vartheta^\lambda(x,\hat{\theta}_N, \hat{r}_N)\right)\right)\\
        &\leq \sum_{j=1}^n{\theta}^*_j\mathbb{E}_{Q_j}[\varphi(x,\xi)]-\lambda\|\theta^*\|^2-\left(\sum_{j=1}^n\tilde{\theta}^*_j\mathbb{E}_{Q_j}[\varphi(x,\xi)] - \lambda\|\tilde{\theta}^*\|^2\right)\\
        &=\nu^\lambda(x,\hat{\theta}_N,\hat{r}_N)-\nu^\lambda(x,\tilde{\theta}_N, \tilde{r}_N)
        \leq (nM+2\lambda)(\|\hat{\theta}_N-\tilde{\theta}_N\|_\infty+\left|\tilde{r}_N-\hat{r}_N\right|).
    \end{aligned}
\end{equation*}
Thus, we have $\mathbb{D}\left(\vartheta^\lambda(x,\tilde{\theta}_N, \tilde{r}_N),\vartheta^\lambda(x,\hat{\theta}_N, \hat{r}_N)\right) \leq (\psi_x^\lambda)^{-1}((nM+2\lambda)(\|\hat{\theta}_N-\tilde{\theta}_N\|_\infty+\left|\tilde{r}_N-\hat{r}_N\right|)).$ The similar result holds for $\mathbb{D}\left(\vartheta^\lambda(x,\hat{\theta}_N, \hat{r}_N),\vartheta^\lambda(x,\tilde{\theta}_N, \tilde{r}_N)\right)$. By the assumed growth condition \(\psi_x^\lambda(\tau)\ge\gamma_\lambda\tau^{p_\lambda}\), we have $(\psi_x^\lambda)^{-1}(t) \leq\left(\frac{t}{\gamma_\lambda}\right)^{1 / p_\lambda}$.  Thus, we can draw the conclusion for $\sup_{x\in\mathcal{X}}\mathbb{H}\left(\vartheta^\lambda(x,\tilde{\theta}_N, \tilde{r}_N),\vartheta^\lambda(x,\hat{\theta}_N, \hat{r}_N)\right)$.

For the second part, let $\tilde{\theta}^*\in\vartheta^\lambda(x,\tilde{\theta}_N, \tilde{r}_N)$ and $\theta^*\in\vartheta(x,\hat{\theta}_N, \hat{r}_N)$ be such that $\|\theta^*-\tilde{\theta}^*\|_\infty=\mathbf{d}\left(\tilde{\theta}^*, \vartheta(x,\hat{\theta}_N, \hat{r}_N)\right)$. This leads to \begin{equation*}
    \begin{aligned}
        \psi_x(\|\theta^*-\tilde{\theta}^*\|_\infty)&=\psi_x\left(\mathbf{d}\left(\tilde{\theta}^*, \vartheta(x,\hat{\theta}_N, \hat{r}_N)\right)\right)\\
        &\leq\nu(x,\hat{\theta}_N,\hat{r}_N)-\nu^\lambda(x,\hat{\theta}_N,\hat{r}_N)+\nu^\lambda(x,\hat{\theta}_N,\hat{r}_N)-\nu^\lambda(x,\tilde{\theta}_N, \tilde{r}_N)\\
        &
        \leq \lambda+(nM+2\lambda)(\|\hat{\theta}_N-\tilde{\theta}_N\|_\infty+\left|\tilde{r}_N-\hat{r}_N\right|).
    \end{aligned}
\end{equation*}
Thus we can obtain the conclusion for $\sup_{x\in\mathcal{X}}\mathbb{D}\left(\vartheta^\lambda(x,\tilde{\theta}_N, \tilde{r}_N),\vartheta(x,\hat{\theta}_N, \hat{r}_N)\right)$.
$\hfill\blacksquare$

To further quantify the perturbation effects on the $x$-solution set, we define a growth function for \eqref{bdrvi-regular2} inspired by \cite{jiang2024distributionally}. For fixed $(\hat{\theta}_N, \hat{r}_N)$, define $\psi^\lambda_{\hat{\theta}_N,\hat{r}_N}: \mathbb{R}_{+} \to \mathbb{R}_{+}$ as
$$
\psi^\lambda_{\hat{\theta}_N,\hat{r}_N}(\tau):=\inf_{x\in\mathcal{X}} \left\{\mathbf{d}\left(0, \sum_{j=1}^n\theta_j\mathbb{E}_{Q_j}[\Phi(x,\xi)]+\mathcal{N}_\mathcal{X}(x)\right):  \theta\in\vartheta^\lambda(x, \hat{\theta}_N,\hat{r}_N), \mathbf{d}\left(x, \mathfrak{X}^\lambda(\hat{\theta}_N,\hat{r}_N)\right) \geq \tau\right\}.
$$
Then its inverse $(\psi_{\hat{\theta}_N,\hat{r}_N}^\lambda)^{-1}$
can be defined similarly as that of $\psi_x^\lambda$.
To ease notation, we use
$\psi_{\hat{\theta}_N,\hat{r}_N}:=\psi^\lambda_{\hat{\theta}_N,\hat{r}_N}$
when $\lambda=0$.

\begin{lemma}\label{lem-con-growth}
For any $\epsilon>0$, $t\geq0$ and fixed $(\hat{\theta}_N,\hat{r}_N)\in\Theta\times[0,1]$, there exists $\delta>0$ such that $(\psi_{\hat{\theta}_N,\hat{r}_N}^\lambda)^{-1}(t+\delta)\leq(\psi_{\hat{\theta}_N,\hat{r}_N}^\lambda)^{-1}(t) + \epsilon$.
\end{lemma}
 \textit{Proof}.
Suppose for contradiction that the conclusion does not hold. Then, there exist an $\epsilon_0>0$ and a sequence $\{\delta_k\}\downarrow0$ such that $(\psi_{\hat{\theta}_N,\hat{r}_N}^\lambda)^{-1}\left(t+\delta_k\right) \geq (\psi_{\hat{\theta}_N,\hat{r}_N}^\lambda)^{-1}(t)+\epsilon_0$ for $k \in \mathbb{N}$.
Specifically, we have
$$
(\psi_{\hat{\theta}_N,\hat{r}_N}^\lambda)^{-1}\left(t+\delta_k\right)=\sup \left\{\tau: \psi_{\hat{\theta}_N,\hat{r}_N}^\lambda(\tau) \leq t+\delta_k\right\} \geq (\psi_{\hat{\theta}_N,\hat{r}_N}^\lambda)^{-1}(t)+\epsilon_0.
$$
Then,
we have that for every $k \in \mathbb{N}$, $\psi^\lambda_{\hat{\theta}_N,\hat{r}_N}\left((\psi_{\hat{\theta}_N,\hat{r}_N}^\lambda)^{-1}(t)+\frac{\epsilon_0}{2}\right) \leq t+\delta_k$. Therefore, as $k \to \infty$, we obtain $\psi^\lambda_{\hat{\theta}_N,\hat{r}_N}\left((\psi_{\hat{\theta}_N,\hat{r}_N}^\lambda)^{-1}(t)+\frac{\epsilon_0}{2}\right) \leq t$, which leads to a contradiction.
$\hfill\blacksquare$

\begin{theorem}[Quantitative stability of $x$-solution set]\label{thm-q-os}
    Suppose that the conditions in Proposition \ref{prop-H} hold. Then, for any $\lambda\geq0$, $\hat{\theta}_N$, $\tilde{\theta}_N \in \Theta$ and $\hat{r}_N, \tilde{r}_N \in [0,1]$, we have
\begin{equation}\label{eq-qs-1}
\mathbb{D}\left( \mathfrak{X}^\lambda(\tilde{\theta}_N,\tilde{r}_N),\mathfrak{X}^\lambda(\hat{\theta}_N,\hat{r}_N)\right)  \leq (\psi_{\hat{\theta}_N,\hat{r}_N}^\lambda)^{-1}\left(\left(\frac{(nM)^{p_\lambda}(nM+2\lambda)}{\gamma_\lambda}(\|\hat{\theta}_N-\tilde{\theta}_N\|_\infty+\left|\tilde{r}_N-\hat{r}_N\right|)\right)^{1/p_\lambda}\right),
\end{equation}
\begin{equation}\label{eq-qs-2}
\mathbb{D}\left( \mathfrak{X}^\lambda(\tilde{\theta}_N,\tilde{r}_N),\mathfrak{X}(\hat{\theta}_N,\hat{r}_N)\right)  \leq (\psi_{\hat{\theta}_N,\hat{r}_N})^{-1}\left(\left(\frac{(nM)^{p_0}(nM+2\lambda)}{\gamma_0}(\|\hat{\theta}_N-\tilde{\theta}_N\|_\infty+\left|\tilde{r}_N-\hat{r}_N\right|)+\frac{\lambda}{\gamma_0}\right)^{1/p_0}\right).
\end{equation}
\end{theorem}
 \textit{Proof}.
For the first part, given any $x \in \mathfrak{X}^\lambda(\tilde{\theta}_N,\tilde{r}_N)$, there exists $\tilde{\theta}^*\in\vartheta^\lambda(x,\tilde{\theta}_N, \tilde{r}_N)$ such that
$$
0 \in \sum_{j=1}^n\tilde{\theta}^*_j\mathbb{E}_{Q_j}[\Phi(x,\xi)]+\mathcal{N}_\mathcal{X}(x).
$$

For any $\hat{\theta}^* \in \vartheta^\lambda(x,\hat{\theta}_N,\hat{r}_N)$, we have
$$
\begin{aligned}
& \left\|\sum_{j=1}^n\tilde{\theta}^*_j\mathbb{E}_{Q_j}[\Phi(x,\xi)]-\sum_{j=1}^n\hat{\theta}^*_j\mathbb{E}_{Q_j}[\Phi(x,\xi)]\right\|_\infty \\
\geq & \mathbf{d}\left(0, \sum_{j=1}^n\hat{\theta}^*_j\mathbb{E}_{Q_j}[\Phi(x,\xi)]+\mathcal{N}_\mathcal{X}(x)\right)-\mathbf{d}\left(0, \sum_{j=1}^n\tilde{\theta}^*_j\mathbb{E}_{Q_j}[\Phi(x,\xi)]+\mathcal{N}_\mathcal{X}(x)\right) \\
= & \mathbf{d}\left(0,\sum_{j=1}^n\hat{\theta}^*_j\mathbb{E}_{Q_j}[\Phi(x,\xi)]+\mathcal{N}_\mathcal{X}(x)\right)
\geq \psi^\lambda_{\hat{\theta}_N,\hat{r}_N}\left(\mathbf{d}\left(x, \mathfrak{X}^\lambda(\hat{\theta}_N,\hat{r}_N)\right)\right),
\end{aligned}
$$
which implies that
$$
\mathbf{d}\left(x, \mathfrak{X}^\lambda(\hat{\theta}_N,\hat{r}_N)\right)\leq (\psi_{\hat{\theta}_N,\hat{r}_N}^\lambda)^{-1}\left( \inf _{\hat{\theta}^* \in \vartheta^\lambda(x,\hat{\theta}_N,\hat{r}_N)}\left\|\sum_{j=1}^n\tilde{\theta}^*_j\mathbb{E}_{Q_j}[\Phi(x,\xi)]-\sum_{j=1}^n\hat{\theta}^*_j\mathbb{E}_{Q_j}[\Phi(x,\xi)]\right\|_\infty\right).
$$

Note that
$$
\begin{aligned}
&\inf_{\hat{\theta}^* \in \vartheta^\lambda(x,\hat{\theta}_N,\hat{r}_N)}\left\|\sum_{j=1}^n\tilde{\theta}^*_j\mathbb{E}_{Q_j}[\Phi(x,\xi)]-\sum_{j=1}^n\hat{\theta}^*_j\mathbb{E}_{Q_j}[\Phi(x,\xi)]\right\|_\infty\\
& \leq nM\mathbf{d}\left(\tilde{\theta}^*, \vartheta^\lambda(x,\hat{\theta}_N,\hat{r}_N)\right) \\
& \leq nM \mathbb{D}\left(\vartheta^\lambda(x,\tilde{\theta}_N,\tilde{r}_N), \vartheta^\lambda(x,\hat{\theta}_N,\hat{r}_N)\right) \\
& \leq \left(\frac{(nM)^{p_\lambda}(nM+2\lambda)}{\gamma_\lambda}(\|\hat{\theta}_N-\tilde{\theta}_N\|_\infty+\left|\tilde{r}_N-\hat{r}_N\right|)\right)^{1/p_\lambda},
\end{aligned}
$$
we obtain $$
\mathbb{D}\left( \mathfrak{X}^\lambda(\tilde{\theta}_N,\tilde{r}_N),\mathfrak{X}^\lambda(\hat{\theta}_N,\hat{r}_N)\right) \leq (\psi_{\hat{\theta}_N,\hat{r}_N}^\lambda)^{-1}\left(\left(\frac{(nM)^{p_\lambda}(nM+2\lambda)}{\gamma_\lambda}(\|\hat{\theta}_N-\tilde{\theta}_N\|_\infty+\left|\tilde{r}_N-\hat{r}_N\right|)\right)^{1/p_\lambda}\right).
$$

For the second part, based on Proposition~\ref{prop-H}, the proof of \eqref{eq-qs-2} follows similarly to that of \eqref{eq-qs-1}, and we omit the details here.
$\hfill\blacksquare$

For any fixed $\lambda\geq0$, applying \eqref{eq-qs-1} twice--once with $(\hat{\theta}_N,\hat{r}_N)$ versus $(\tilde{\theta}_N,\tilde{r}_N)$, and once with the roles swapped--immediately yields the Hausdorff bound
\begin{equation}\label{eq-haus-x}
    \mathbb{H}\left( \mathfrak{X}^\lambda(\tilde{\theta}_N,\tilde{r}_N),\mathfrak{X}^\lambda(\hat{\theta}_N,\hat{r}_N)\right)  \leq \max\left\{(\psi_{\hat{\theta}_N,\hat{r}_N}^\lambda)^{-1}((C_\lambda\Delta_N)^{\frac{1}{p_\lambda}}),(\psi_{\tilde{\theta}_N,\tilde{r}_N}^\lambda)^{-1}((C_\lambda\Delta_N)^{\frac{1}{p_\lambda}})\right\},
\end{equation}
where $\Delta_N:=\|\hat{\theta}_N-\tilde{\theta}_N\|_\infty+\left|\tilde{r}_N-\hat{r}_N\right|$ and $C_\lambda:=\frac{(nM)^{p_\lambda}(nM+2\lambda)}{\gamma_\lambda}$.
We next derive two immediate consequences of Theorem~\ref{thm-q-os}.
The first provides a posterior finite-sample guarantee for the regularized data-driven solution set, {\color{black}while the second gives a direct deterministic bound for the perturbed regularized solution set relative to the true Bayesian DRVI solution set.}

\begin{corollary}[Posterior finite-sample guarantee]\label{coro-1}
Suppose that the conditions in Proposition \ref{prop-H} hold. For any given finite sample set $\mathcal{S}^N$,  with $C_0:=\frac{(nM)^{p_0}(nM+2\lambda)}{\gamma_0}$, the following guarantee holds:
$$
\mathbb{P}_{\mathcal{S}^N}\left\{\mathbb{D}\left( \mathfrak{X}^\lambda(\hat{\theta}_N,\hat{r}_N),\mathfrak{X}(\theta^c,r^c)\right) \leq (\psi_{\theta^c,r^c})^{-1}\left(\left(C_0(\hat{r}_N+\hat{\delta}_N)+\frac{\lambda}{\gamma_0}\right)^{1/p_0}\right)
\right\}\geq1-\alpha.
$$
\end{corollary}
 \textit{Proof}.
For any $x \in \mathfrak{X}^\lambda(\hat{\theta}_N,\hat{r}_N)$, there exists $\hat{\theta}^*\in\vartheta^\lambda(x,\hat{\theta}_N, \hat{r}_N)$ such that
$
0 \in \sum_{j=1}^n\hat{\theta}^*_j\mathbb{E}_{Q_j}[\Phi(x,\xi)]+\mathcal{N}_\mathcal{X}(x).
$
Then, for any $\theta^*\in\vartheta(x,\theta^c,r^c)$, by the same argument as in the proof of Theorem~\ref{thm-q-os},
$$
\mathbf{d}\left(x, \mathfrak{X}(\theta^c,r^c)\right)\leq (\psi_{\theta^c,r^c})^{-1}\left( \left\|\sum_{j=1}^n\hat{\theta}^*_j\mathbb{E}_{Q_j}[\Phi(x,\xi)]-\sum_{j=1}^n{\theta}^*_j\mathbb{E}_{Q_j}[\Phi(x,\xi)]\right\|_\infty \right).
$$
In the event that $\Theta^c\subseteq\hat{\Theta}_N$, we have
$\|\theta^c-\hat{\theta}_N\|_\infty\leq\hat{r}_N$ and $|r^c-\hat{r}_N|=\hat{\delta}_N$.
Note that
$$
\left\|\sum_{j=1}^n\hat{\theta}^*_j\mathbb{E}_{Q_j}[\Phi(x,\xi)]-\sum_{j=1}^n{\theta}^*_j\mathbb{E}_{Q_j}[\Phi(x,\xi)]\right\|_\infty \leq \left(C_0(\hat{r}_N+\hat{\delta}_N)+\frac{\lambda}{\gamma_0}\right)^{1/p_0},
$$
we obtain $$
\begin{aligned}
    \mathbb{D}\left( \mathfrak{X}^\lambda(\hat{\theta}_N,\hat{r}_N),\mathfrak{X}(\theta^c,r^c)\right)
    &\leq (\psi_{\theta^c,r^c})^{-1}\left(\left(C_0(\hat{r}_N+\hat{\delta}_N)+\frac{\lambda}{\gamma_0}\right)^{1/p_0}\right).
\end{aligned}
$$
In terms of probability, based on Theorem \ref{thm-new}, we have
$$
\begin{aligned}
    \mathbb{P}_{\mathcal{S}^N}\left\{\mathbb{D}\left( \mathfrak{X}^\lambda(\hat{\theta}_N,\hat{r}_N),\mathfrak{X}(\theta^c,r^c)\right) \leq (\psi_{\theta^c,r^c})^{-1}\left(\left(C_0(\hat{r}_N+\hat{\delta}_N)+\frac{\lambda}{\gamma_0}\right)^{1/p_0}\right)
    \right\}\geq \mathbb{P}_{\mathcal{S}^N}\left\{\Theta^c\subseteq\hat{\Theta}_N\right\}\geq1-\alpha,
\end{aligned}
$$
which completes the proof.
$\hfill\blacksquare$

{
\color{black}
Corollary \ref{coro-1} gives a posterior finite-sample guarantee based on the Bayesian ambiguity-set construction.
In addition, by replacing $(\hat\theta_N,\hat r_N)$ to $(\theta^c,r^c)$ in \eqref{eq-qs-2}, we obtain the following direct deterministic bound to the true solution set.
\begin{corollary}\label{cor-direct-truth}
Suppose that the conditions in Proposition \ref{prop-H} hold. Then, for any contaminated sample set $\tilde{\mathcal S}^N$ with induced ambiguity parameters $(\tilde\theta_N,\tilde r_N)$,
\[
\mathbb{D}\left( \mathfrak{X}^\lambda(\tilde{\theta}_N,\tilde{r}_N),\mathfrak{X}(\theta^c,r^c)\right)
\leq
(\psi_{\theta^c,r^c})^{-1}\left(
\left(
C_0
\bigl(\|\tilde{\theta}_N-\theta^c\|_\infty+|\tilde r_N-r^c|\bigr)
+\frac{\lambda}{\gamma_0}
\right)^{1/p_0}
\right).
\]
\end{corollary}

}

{\color{black}
The following remark complements the sequential $\lim_{N\to\infty}\lim_{\lambda\downarrow0}$ in Theorem~\ref{thm-con}
by providing a data-dependent rule for coupling $\lambda_N$ with the sample size $N$.

\begin{remark}
By applying \eqref{eq-qs-2} with $(\theta^c,r^c)$ and
$(\hat\theta_N,\hat r_N)$, we obtain
\[
\mathbb D\left(\mathfrak X^{\lambda_N}(\hat\theta_N,\hat r_N),\mathfrak X(\theta^c,r^c)\right)
\le
(\psi_{\theta^c,r^c})^{-1}
\left(\left(\frac{(nM)^{p_0}(nM+2\lambda_N)}{\gamma_0}(\|\hat{\theta}_N-\theta^c\|_\infty+\hat{\delta}_N)+\frac{\lambda_N}{\gamma_0}\right)^{1/p_0}\right).
\]
This bound separates the total approximation error into two components:
the estimation error
$\|\hat{\theta}_N-\theta^c\|_\infty+\hat{\delta}_N$,
which is induced by finite-sample uncertainty in the ambiguity-set construction,
and the regularization bias
$\lambda_N$.
Hence the bound yields a bias-variance type trade-off between statistical estimation error in the ambiguity-set construction and regularization bias in the lower-level problem, in the spirit of standard regularization analyses; see, e.g., \cite[Section~7.6]{bach2024learning} for an analogous discussion.
In particular, under the construction in Section~2 with $\|\hat{\theta}_N-\theta^c\|_\infty\to0$ and $\hat{\delta}_N\to0$, any choice satisfying $\lambda_N\to0$ yields consistency of the regularized data-driven solution set. Moreover,  under the Bernstein--von Mises approximation, one typically has $\|\hat{\theta}_N-\theta^c\|_\infty=O(N^{-1/2})$ and $\hat \delta_N=O(N^{-1/2})$, then choosing
$\lambda_N\asymp N^{-1/2}$ balances the sampling error and the regularization bias at the same order.
\end{remark}
}

In practice, data often deviate from the assumed true distribution $P_{\theta^c}$ due to contamination or outliers, especially in fields such as finance, biology and medicine \cite{lecue2020robust}.
Typically, one observes samples $\tilde{\mathcal{S}}^N = \{ \tilde{\xi}^1, \ldots, \tilde{\xi}^N \}$ generated from a perturbed distribution $P_{\tilde{\theta}^c}$ rather than $P_{\theta^c}$.
Following \cite{pichler2022quantitative}, we assume that $\tilde{\xi}^1,\dots,\tilde{\xi}^N$ are i.i.d. from $P_{\tilde{\theta}^c}$. In what follows, we further examine how the discrepancy between $P_{\tilde{\theta}^c}$ and $P_{\theta^c}$ affects the $x$-solution set of the Bayesian DRVI problem.
{\color{black}
In the following, whenever the perturbed sample sequence
$\tilde{\mathcal S}^N$ is considered, the corresponding
``almost surely'' statements are understood with respect to
$(P_{\tilde\theta^c})^{\mathbb N}$. Statements involving both the original and perturbed sample sequences are understood
almost surely on the joint probability space generated by these two sequences.}

\begin{lemma}\label{lem-4}
Suppose that Assumptions \ref{ass-bvm}-\ref{ass3.1} and the conditions in Lemma \ref{lem-lip-optimal} hold.
Moreover, for any fixed $\lambda\geq0$, there exist a neighborhood $U$ of $(\theta^c,r^c)$ in $\Theta\times[0,1]$, constants $\bar p_\lambda\geq1$, $\bar\gamma_\lambda>0$, and $\tau_0>0$ such that for all $(\theta,r)\in U$ and $\tau\in[0,\tau_0]$,
\[
\psi_{\theta,r}^\lambda(\tau)\ge\bar\gamma_\lambda\tau^{\bar p_\lambda}.
\]
Then, for any  $\tau > 0$ and $t > 0$, we have (i) $\liminf_{N \to \infty} \psi^\lambda_{\hat{\theta}_N, \hat{r}_N}(\tau) \geq \psi_{\theta^c,r^c}^\lambda(\tau)$, and (ii) $\limsup_{N \to \infty}( \psi_{\hat{\theta}_N, \hat{r}_N}^\lambda)^{-1}(t) \leq (\psi^\lambda_{\theta^c,r^c})^{-1}(t)$.
\end{lemma}

 \textit{Proof}. We prove part (i) by contradiction. Suppose that there exists $\delta > 0$ such that \[\liminf_{N\to\infty}\psi^\lambda_{\hat{\theta}_N,\hat{r}_N}(\tau)\leq\psi_{\theta^c,r^c}^\lambda(\tau)-\delta.\] Thus, there exists a subsequence $\left\{N_k\right\}_{k \geq 1}$ such that $\psi^{\lambda}_{\hat{\theta}_{N_k},\hat{r}_{N_k}}(\tau)\leq\psi_{\theta^c,r^c}^\lambda(\tau)-\delta.$ By the definition of $\psi^{\lambda}_{\hat{\theta}_{N_k},\hat{r}_{N_k}}(\tau)$, there exists a pair $(x_{N_k},{\theta}^*_{N_k})$ such that $\mathbf{d}\left(x_{N_k}, \mathfrak{X}^{\lambda}(\hat{\theta}_{N_k},\hat{r}_{N_k})\right) \geq \tau$ and \begin{equation}\label{eq-123}
 \mathbf{d}\left(0, \sum_{j=1}^n{\theta}^*_{N_k,j}\mathbb{E}_{Q_j}[\Phi(x_{N_k},\xi)]+\mathcal{N}_\mathcal{X}(x_{N_k})\right)\leq\psi_{\theta^c,r^c}^\lambda(\tau)-\delta.
 \end{equation}
 Letting $x_{N_k} \to x^*$, we claim $\mathbf{d}(x^*, \mathfrak{X}^\lambda(\theta^c, r^c)) \geq \tau$. Otherwise, there exists $x^\dagger \in \mathfrak{X}^{\lambda}(\theta^c, r^c)$ such that $\|x^* - x^\dagger\| < \tau$. From \eqref{eq-haus-x}, we have $$\mathbf{d}\left( x^\dagger,\mathfrak{X}^\lambda(\hat{\theta}_N,\hat{r}_N)\right)\leq( \psi_{\hat{\theta}_N, \hat{r}_N}^\lambda)^{-1}\left(C_\lambda^{1/p_\lambda}(\|\hat{\theta}_N-\theta^c\|_\infty+|\hat{r}_N-r^c|)^{1/p_\lambda}\right):=\epsilon_N.$$
 From Theorem \ref{thm-new} the assumed growth condition, for $N$ large enough we have $(\hat{\theta}_N,\hat r_N)\in U$ and $
\epsilon_N\to\ 0$ as $N\to\infty$.
For each $N_k$ large enough, there exists $y_{N_k}\in \mathfrak{X}^{\lambda}\left(\hat{\theta}_{N_k}, \hat{r}_{N_k}\right)$ with $\left\|y_{N_k}-x^{\dag}\right\|_{\infty}<\epsilon_{N_k}. $ Thus, we have
 $$
 \begin{aligned} \left\|x_{N_k}-y_{N_k}\right\|_{\infty} \leq\left\|x_{N_k}-x^*\right\|_{\infty}+\left\|x^*-x^{\dagger}\right\|_{\infty}+\left\|x^{\dagger}-y_{N_k}\right\|_{\infty}<\left(\epsilon^{\prime}+\epsilon_{N_k}\right)+\left\|x^*-x^{\dagger}\right\|_{\infty}
 \end{aligned}
 $$
 where $\epsilon^{\prime} \to 0$ as $N_k \to \infty$. Since $\left\|x^*-x^{\dagger}\right\|_{\infty}<\tau$ and $\epsilon_{N_k}$ can be made arbitrarily small, $\left\|x_{N_k}-y_{N_k}\right\|_{\infty}$ can be made less than $\tau$ for large enough $N_k$. This contradicts the assumption that $\mathbf{d}\left(x_{N_k}, \mathfrak{X}^{\lambda}\left(\hat{\theta}_{N_k}, \hat{r}_{N_k}\right)\right) \geq \tau$. Therefore, our assumption must be false, implying $\mathbf{d}\left(x^*, \mathfrak{X}^\lambda\left(\theta^c,r^c\right)\right) \geq \tau$. Using an argument similar to that in Lemma \ref{lem-con-svi}, it is known from \eqref{eq-123} that $$
 \mathbf{d}\left(0, \sum_{j=1}^n \theta^*_j \mathbb{E}_{Q_j}[\Phi\left(x^*, \xi\right)]+\mathcal{N}_{\mathcal{X}}\left(x^*\right)\right) \leq \psi_{\theta^c,r^c}^\lambda(\tau)-\frac{\delta}{2}. $$
 However, by the definition of $\psi_{\theta^c,r^c}^\lambda(\tau)$, this is impossible because $\psi_{\theta^c,r^c}^\lambda(\tau)$ is the infimum of such distances for points at least $\tau$ away from $\mathfrak{X}^{\lambda}\left(\theta^c, r^c\right)$. Thus, $$ \psi_{\theta^c,r^c}^\lambda(\tau) \leq \mathbf{d}\left(0, \sum_{j=1}^n \theta^*_j \mathbb{E}_{Q_j}[\Phi\left(x^*, \xi\right)]+\mathcal{N}_{\mathcal{X}}\left(x^*\right)\right) \leq \psi_{\theta^c,r^c}^\lambda(\tau)-\frac{\delta}{2}, $$ which is a contradiction. Therefore, the inequality $ \liminf _{N \to \infty} \psi^\lambda_{\hat{\theta}_N, \hat{r}_N}(\tau) \geq \psi_{\theta^c,r^c}^\lambda(\tau)$ must hold.

To establish part (ii), note that for any $\epsilon > 0$, part (i) implies the existence of $N_0$ such that for all \(N \geq N_0\),
$
\psi^\lambda_{\hat{\theta}_N, \hat{r}_N}(\tau) \geq \psi_{\theta^c,r^c}^\lambda(\tau) - \epsilon.
$
By the definition of $(\psi^\lambda_{\hat{\theta}_N, \hat{r}_N})^{-1}$ and $(\psi_{\theta^c,r^c}^\lambda)^{-1}$, we have that
$ (\psi_{\theta^c,r^c}^\lambda)^{-1}(t+\epsilon) \geq (\psi^\lambda_{\hat{\theta}_N, \hat{r}_N})^{-1}(t).$
Since $\epsilon$ is arbitrary, we conclude that
$
\limsup_{N \to \infty}( \psi_{\hat{\theta}_N, \hat{r}_N}^\lambda)^{-1}(t)  \leq (\psi_{\theta^c,r^c}^\lambda)^{-1}(t),
$
which completes the proof. $\hfill\blacksquare$

{\color{black}

}
Recall that Theorem \ref{thm-q-os} establishes quantitative stability for the Bayesian DRVI problem. In conjunction with Lemma \ref{lem-con-growth}, it offers explicit upper bounds for the qualitative results derived in Theorem \ref{thm-con}. To maintain self-containment, we will now demonstrate that, for any $\lambda\geq0$, the right-hand side of equation \eqref{eq-qs-1} tends to 0 as $N\to\infty$ when the perturbation $\|\theta^c-\tilde{\theta}^c\|_\infty$ tends to 0.
\begin{theorem}\label{thm-10}
    Suppose that the conditions in Lemma \ref{lem-4} hold. Then, for any $\lambda\geq0$, almost surely
    $$
\lim_{N\to\infty,\|{\theta}^c-\tilde{\theta}^c\|_\infty\to0}\mathbb{D}\left( \mathfrak{X}^\lambda(\tilde{\theta}_N,\tilde{r}_N),{\mathfrak{X}^\lambda(\theta^c,r^c)}\right) =0.
$$
\end{theorem}
 \textit{Proof}.
   Note that $$
   \begin{aligned}
&\limsup_{N\to\infty,\|{\theta}^c-\tilde{\theta}^c\|_\infty\to0}\mathbb{D}\left( \mathfrak{X}^\lambda(\tilde{\theta}_N,\tilde{r}_N),{\mathfrak{X}^\lambda(\theta^c,r^c)}\right) \\
\leq
&\limsup_{N\to\infty,\|{\theta}^c-\tilde{\theta}^c\|_\infty\to0}(\psi_{\theta^c,r^c}^\lambda)^{-1}\left(\left(C_\lambda(\|\tilde{\theta}_N-{\theta}^c\|_\infty+\left|\tilde{r}_N-r^c\right|)\right)^{1/p_\lambda}\right)\\
\leq&\limsup_{\|{\theta}^c-\tilde{\theta}^c\|_\infty\to0}(\psi_{\theta^c,r^c}^\lambda)^{-1}\left(\left(C_\lambda\|\theta^c-\tilde{\theta}^c\|_\infty\right)^{1/p_\lambda}\right)=0.
   \end{aligned}
$$
Since $(\psi_{\theta^c,r^c}^\lambda)^{-1}\left(\left(C_\lambda\|\theta^c-\tilde{\theta}^c\|_\infty\right)^{1/p_\lambda}\right)\geq0$, we complete the proof.
$\hfill\blacksquare$

Theorem \ref{thm-10} demonstrates the asymptotic convergence of the $x$-solution set for the contaminated problem \eqref{p-bdrvi2} as the distance between the contaminated distribution $P_{\tilde{\theta}^c}$ and the true distribution $P_{{\theta}^c}$ approaches zero. We now derive a quantitative error bound for the deviation between $x$-solution sets when the  distance between $\theta^c$ and $\tilde{\theta}^c$ is not zero. To this end, we need the following lemma.

\begin{lemma}\label{lem-9}
        Suppose that the conditions in Theorem \ref{thm-10} hold and $\psi_{\theta^c,r^c}^\lambda$ is continuous for any fixed $\lambda\geq0$. Then, for any $\tau > 0$, we have (i) $\limsup_{N\to\infty}\psi_{\hat{\theta}_N,\hat{r}_N}^\lambda(\tau)\leq\psi_{\theta^c,r^c}^\lambda(\tau)$ and (ii) $\liminf_{N\to\infty} (\psi_{\hat{\theta}_N, \hat{r}_N}^\lambda)^{-1}(t) \geq (\psi_{\theta^c,r^c}^\lambda)^{-1}(t)$.
\end{lemma}
 \textit{Proof}.
To prove (i), let $\tau>\delta>0$ and $\epsilon>0$ be given.  By the definition of $\psi_{\theta^c,r^c}^\lambda$, there exist $x^* \in \mathcal{X}$ and $\theta^*\in\vartheta^\lambda\left(x^*, \theta^c, r^c\right)$ such that $\mathbf{d}\left(x^*, \mathfrak{X}^\lambda\left(\theta^c,r^c\right)\right) \geq \tau+\delta$ and $$\mathbf{d}\left(0, \sum_{j=1}^n \theta_j^* \mathbb{E}_{Q_j}[\Phi\left(x^*, \xi\right)]+\mathcal{N}_{\mathcal{X}}\left(x^*\right)\right) \leq \psi_{\theta^c,r^c}^\lambda(\tau+\delta)+\epsilon.$$
Since $(\hat{\theta}_N,\hat{r}_N)\to(\theta^c,r^c)$ and $\mathfrak{X}^\lambda$ is upper semicontinuous, there exists $N_1$ such that for all
$N>N_1$, we have $\mathfrak{X}^\lambda(\hat{\theta}_N,\hat{r}_N)\subseteq\mathfrak{X}^\lambda(\theta^c,r^c)+\frac{\delta}{2}\mathbb{B}_\infty.$
Hence, for any $x'\in\mathfrak{X}^\lambda(\hat{\theta}_N,\hat{r}_N)$, there exists $y\in{\mathfrak{X}^\lambda(\theta^c,r^c)}$ such that $\|x'-y\|_\infty\leq\frac{\delta}{2}$, which implies $$\|x^*-x'
\|_\infty\geq\|x^*-y
\|_\infty-\|y-x'
\|_\infty\geq\tau+\frac{\delta}{2}.$$
Thus, $\mathbf{d}\left(x^*, \mathfrak{X}^\lambda\left(\hat{\theta}_N,\hat{r}_N\right)\right) \geq \tau.$

On the other hand, due to the upper semicontinuity of $\vartheta^\lambda(x,\cdot,\cdot)$ from Theorem \ref{thm-uppercon-1} and Proposition \ref{prop-H}, there exist $N_2\in\mathbb{N}$ and a sequence $\{\theta_N^*\}_{N\geq N_2}$ with $\theta_N^*\in \vartheta^\lambda(x^*, \hat{\theta}_N, \hat{r}_N)$ such that ${\theta}_N^* \to \theta^*$ as $N\to\infty$. Then, we have
$
\psi_{\hat{\theta}_N, \hat{r}_N}^\lambda(\tau) \leq \mathbf{d}\left(0, \sum_{j=1}^n {\theta}_{N,j}^* \mathbb{E}_{Q_j}[\Phi(x^*, \xi)] + \mathcal{N}_\mathcal{X}(x^*)\right).
$
Taking the limit yields
$$
\begin{aligned}
    \limsup _{N \to \infty} \psi^\lambda_{\hat{\theta}_N, \hat{r}_N}(\tau)
    &\leq \lim_{N \to \infty} \mathbf{d}\left(0, \sum_{j=1}^n {\theta}_{N,j}^* \mathbb{E}_{Q_j}[\Phi\left(x^*, \xi\right)]+\mathcal{N}_{\mathcal{X}}\left(x^*\right)\right)\\
    &=\mathbf{d}\left(0, \sum_{j=1}^n \theta_j^* \mathbb{E}_{Q_j}[\Phi\left(x^*, \xi\right)]+\mathcal{N}_{\mathcal{X}}\left(x^*\right)\right)\leq \psi_{\theta^c,r^c}^\lambda(\tau+\delta)+\epsilon.
\end{aligned}
$$
Since $\epsilon>0$ and $\delta>0$ are arbitrary, we conclude
$
\limsup _{N \to \infty} \psi^\lambda_{\hat{\theta}_N, \hat{r}_N}(\tau) \leq \psi_{\theta^c,r^c}^\lambda(\tau).
$

To prove (ii), let \(t \geq 0\) be given and define
$\tau^* = (\psi_{\theta^c,r^c}^\lambda)^{-1}(t) = \sup \left\{ \tau \geq 0 | \psi_{\theta^c,r^c}^\lambda(\tau) \leq t \right\}.$
Since \(\psi_{\theta^c,r^c}^\lambda\) is continuous and non-decreasing, we have
$\psi_{\theta^c,r^c}^\lambda(\tau^*) = t$ and
$\psi_{\theta^c,r^c}^\lambda(\tau) < t$
for any \(\tau < \tau^*\).
By the result of part (i), \(\limsup_{N\to\infty} \psi_{\hat{\theta}_N, \hat{r}_N}^\lambda(\tau) \leq \psi_{\theta^c,r^c}^\lambda(\tau)\). This means that for the fixed \(\tau < \tau^*\) and any \(\epsilon > 0\), there exists $N_0\in\mathbb N$ such that for all \(N \geq N_0\), we have
$
\psi_{\hat{\theta}_N, \hat{r}_N}^\lambda(\tau) \leq \psi_{\theta^c,r^c}^\lambda(\tau) + \epsilon.
$

Since \(\psi_{\theta^c,r^c}^\lambda(\tau) < t\), we can choose \(\epsilon\) small enough such that $
\psi_{\theta^c,r^c}^\lambda(\tau) + \epsilon < t.
$
Thus, for all \(N \geq N_0\), $\psi_{\hat{\theta}_N, \hat{r}_N}^\lambda(\tau) < t$,
which means that for any \(\tau < \tau^*\) and sufficiently large \(N\),
$
(\psi_{\hat{\theta}_N, \hat{r}_N}^\lambda)^{-1}(t) \geq \tau.
$
Since the above inequality holds for all \(\tau < \tau^*\), we have
$
\liminf_{N\to\infty} (\psi_{\hat{\theta}_N, \hat{r}_N}^\lambda)^{-1}(t) \geq \tau^*=(\psi_{\theta^c,r^c}^\lambda)^{-1}(t).
$
$\hfill\blacksquare$

Based on these conclusions, we can further derive the error bound corresponding to Theorem \ref{thm-10}.
\begin{theorem}\label{thm-last}
Suppose that the conditions in Theorem \ref{thm-10} hold and $\psi_{\theta^c,r^c}^\lambda$ is continuous for any fixed $\lambda\geq0$. Then, almost surely
    $$
\lim_{N\to\infty}\mathbb{D}\left( \mathfrak{X}^\lambda(\tilde{\theta}_N,\tilde{r}_N),\mathfrak{X}^\lambda(\hat{\theta}_N,\hat{r}_N)\right) \leq (\psi^\lambda_{\theta^c,r^c})^{-1}\left(\left(C_\lambda\|{\theta}^c-\tilde{\theta}^c\|_\infty\right)^{1/p_\lambda}\right).
$$
\end{theorem}
 \textit{Proof}.
  By Theorem~\ref{thm-q-os}, for any $N$ and $\lambda\geq 0$,
$$
\mathbb{D}\left( \mathfrak{X}^\lambda(\tilde{\theta}_N,\tilde{r}_N),\mathfrak{X}^\lambda(\hat{\theta}_N,\hat{r}_N)\right)  \leq (\psi_{\hat{\theta}_N,\hat{r}_N}^\lambda)^{-1}\left(\left(C_\lambda\Delta_N\right)^{1/p_\lambda}\right).
$$
Taking the limit on both sides of the inequality, for any $\delta\geq\Delta_N-\|{\theta}^c-\tilde{\theta}^c\|_\infty$,  we have that
$$
\begin{aligned}
  \lim_{N\to\infty}\mathbb{D}\left( \mathfrak{X}^\lambda(\tilde{\theta}_N,\tilde{r}_N),\mathfrak{X}^\lambda(\hat{\theta}_N,\hat{r}_N)\right)
  \leq&\lim_{N\to\infty}(\psi_{\hat{\theta}_N,\hat{r}_N}^\lambda)^{-1}\left(\left(C_\lambda\Delta_N\right)^{1/p_\lambda}\right)\\
\leq&\lim_{N\to\infty}(\psi_{\hat{\theta}_N,\hat{r}_N}^\lambda)^{-1}\left(\left(C_\lambda(\delta+\|{\theta}^c-\tilde{\theta}^c\|_\infty)\right)^{1/p_\lambda}\right)\\
\leq&(\psi_{\theta^c,r^c}^\lambda)^{-1}\left(\left(C_\lambda\|{\theta}^c-\tilde{\theta}^c\|_\infty\right)^{1/p_\lambda}\right),
\end{aligned}
$$
where the second inequality holds due to the non-decreasing property of $(\psi_{\hat{\theta}_N,\hat{r}_N}^\lambda)^{-1}$ and
the last inequality follows from Lemma \ref{lem-9}.
Since $\delta$ can be taken arbitrarily small as $N\to\infty$, we obtain that as $N\to\infty$,
$$
\begin{aligned}
(\psi_{\hat{\theta}_N,\hat{r}_N}^\lambda)^{-1}\left(\left(C_\lambda\Delta_N\right)^{1/p_\lambda}\right)
    \leq(\psi_{\theta^c,r^c}^\lambda)^{-1}\left(\left(C_\lambda\|{\theta}^c-\tilde{\theta}^c\|_\infty\right)^{1/p_\lambda}\right).
\end{aligned}
$$
Thus the proof is completed.
$\hfill\blacksquare$

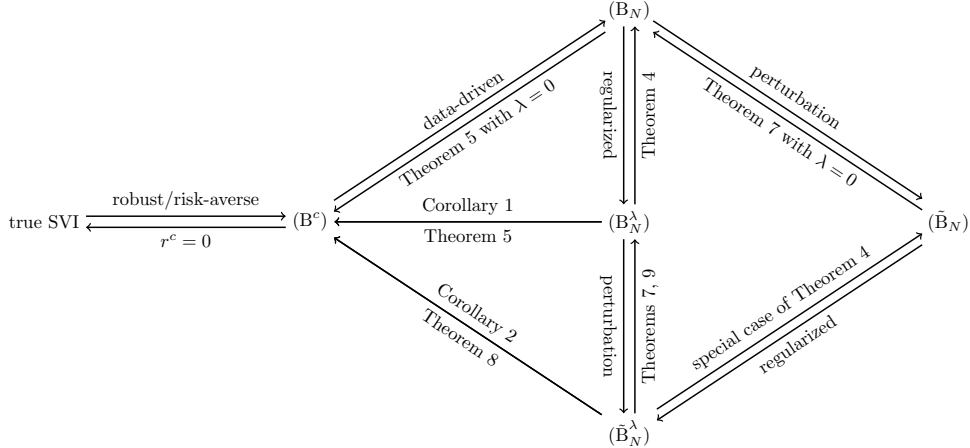
\begin{figure}[h]
\centering
\scalebox{0.7}{
\begin{tikzpicture}
    \node (P_true) at (-10,0) {true SVI};
    \node (P) at (-5,0) {$(\text{B}^c)$};
    \node (Q_alpha) at (7,0) {$(\tilde{\text{B}}_N)$};
    \node (P_alpha_N) at (1,4) {$(\text{B}_N)$};
    \node (r_alpha) at (1,-4) {$(\tilde{\text{B}}_N^\lambda)$};
    \node (r_alpha_N) at (1,0) {$(\text{B}_N^\lambda)$};

    \draw[->, thick, transform canvas={yshift=3pt}] (P_true) -- (P) node[midway, sloped,above] {robust/risk-averse};
    \draw[->, thick, transform canvas={yshift=-3pt}] (P) -- (P_true) node[midway, sloped,below] {$r^c=0$};

    \draw[->, thick] (r_alpha_N)-- (P) node[midway,sloped, below] {Theorem \ref{thm-con}};
    \draw[<-, thick] (P)--(r_alpha_N) node[midway,sloped, above] {Corollary \ref{coro-1}};

    \draw[->, thick] (r_alpha)-- (P) node[midway,sloped, below] {Theorem \ref{thm-10} with $\lambda=0$};

    \draw[->, thick, transform canvas={xshift=-3pt}] (P_alpha_N) -- (r_alpha_N) node[midway, sloped,below] {regularized};
    \draw[->, thick, transform canvas={xshift=3pt}] (r_alpha_N)-- (P_alpha_N) node[midway,sloped, below] {Theorem \ref{thm-con-dr-regular}};

    \draw[->, thick, transform canvas={yshift=3pt}] (P) -- (P_alpha_N) node[midway, sloped,above] {data-driven};
    \draw[->, thick, transform canvas={yshift=-3pt}] (P_alpha_N) -- (P) node[midway,sloped,below] {Theorem \ref{thm-con} with $\lambda=0$};

    \draw[->, thick, transform canvas={yshift=3pt}] (P_alpha_N) -- (Q_alpha) node[midway, sloped, above] {perturbation};
        \draw[->, thick, transform canvas={yshift=-3pt}] (Q_alpha) --(P_alpha_N)  node[midway, sloped, below] {Theorem \ref{thm-q-os} with $\lambda=0$};

    \draw[->, thick, transform canvas={yshift=3pt}] (r_alpha) -- (Q_alpha) node[midway, sloped, above] {special case of Theorem \ref{thm-con-dr-regular}};
        \draw[->, thick, transform canvas={yshift=-3pt}] (Q_alpha) --(r_alpha)  node[midway, sloped, below] {regularized};

            \draw[->, thick, transform canvas={xshift=3pt}] (r_alpha) -- (r_alpha_N) node[midway, sloped, below]{Theorems \ref{thm-q-os}, \ref{thm-last}} ;
        \draw[->, thick, transform canvas={xshift=-3pt}] (r_alpha_N) --(r_alpha)  node[midway, sloped, below]{perturbation} ;
\end{tikzpicture}}
\caption{The relationship diagram of problems \eqref{bdrvi-true}, \eqref{bdrvi2}, \eqref{bdrvi-regular2}, \eqref{p-bdrvi2}.}\label{fig:relationship}
\end{figure}

To clarify the relationships among different Bayesian DRVI problems discussed in Sections 3 and 4, we draw a relationship diagram in Figure~\ref{fig:relationship}.
Here, $(\tilde{\text{B}}_N)$ denotes the perturbed formulation \eqref{p-bdrvi2} with $\lambda=0$.
Overall, Theorem~\ref{thm-con} and Corollary~\ref{coro-1} show that a properly designed Bayesian parametric ambiguity set--formulated in Theorem~\ref{thm-new}--ensures both finite-sample and asymptotic guarantees for the $x$-solution set of the regularized data-driven problem \eqref{bdrvi-regular2}.
In addition, Theorems~\ref{thm-q-os} and \ref{thm-last} establish quantitative stability of the
$x$-solution set between \eqref{bdrvi-regular2} and its perturbed counterpart \eqref{p-bdrvi2} under sample contamination.
Theorem~\ref{thm-10} further demonstrates that the perturbed problem \eqref{p-bdrvi2} asymptotically converges to the Bayesian DRVI formulation \eqref{bdrvi-true}, thereby confirming the robustness of the
$x$-solution set against data contamination.

\section{Numerical experiment}
In this section, we illustrate the effectiveness of the proposed Bayesian DRVI model \eqref{bdrvi-regular2}.
We consider the distributionally robust multi-portfolio Nash equilibrium problem in Example \ref{ex-final} with $I=2$ accounts, described as follows:
\begin{equation}
    \min_{ x^i \in \mathcal{X}_i} \max _{\theta^i \in \hat{\Theta}_{N}} \sum_{j=1}^n\theta_{j}^i\mathbb{E}_{Q_j}\left[\varphi\left(x^i,  \xi\right)\right] +\frac{\kappa}{2}\|x^i-x^{-i}\|^2,\ i=1,  2 . \label{numer-nash}
\end{equation}
Here, $\mathcal{X}_i=\{x^i\mid e^\top x^i\leq1,\ x^i\in\mathbb{R}_+^d\}$ denotes the feasible portfolio set, where the inequality allows investment in a risk-free asset, and $ \varphi(x^i, \xi):=-\log(1+{x^i}^\top\xi)$ represents the logarithmic disutility function which is convex in $x^i$. The term $\frac{\kappa}{2}\|x^i-x^{-i}\|^2$, with $\kappa>0$ being a common competition parameter penalizing deviations between the two portfolios (see \cite{lacker2019mean}).
In the following, let $x=(x^1,x^2)$ and denote
$$
\begin{aligned}
&\varphi^n(x^i):=\left(\mathbb{E}_{Q_1}[\varphi(x^i,\xi)],\ldots,\mathbb{E}_{Q_n}[\varphi(x^i,\xi)]\right)^\top,
   &\  &u_i(x^i,\theta^i):={\theta^i}^\top \varphi^n(x^i)-\lambda\|\theta^i\|^2,\\
   &U_i(x^i):=\max _{\theta^i \in \hat{\Theta}_{N}} u_i(x^i,\theta^i), &\ &f_\lambda(x):=U_1(x^1)+U_2(x^2)
     +\tfrac{\kappa}{2}\|x^1-x^{2}\|^2.
\end{aligned}
$$
As shown in Example \ref{ex-final}, \eqref{numer-nash} leads to the following regularized Bayesian DRVI formulation:
\begin{equation}\label{num-eq-1}
    \begin{aligned}
& 0 \in \sum_{j=1}^n\theta_{j}^i\mathbb{E}_{Q_j}[\nabla_{x^i}\varphi(x^i,\xi)]+\kappa\bigl(x^i-x^{-i}\bigr)+\mathcal{N}_{\mathcal{X}_i}\left(x^i\right),\ i=1,2,\\
& \theta^i\in\underset{{\theta}^{\prime}\in\hat{\Theta}_N}{\argmax}\ u_i(x^i,{\theta}^{\prime}),\ i=1,2.
\end{aligned}
\end{equation}
{\color{black}
For each fixed $x^i\in\mathcal X_i$, the lower-level objective
$u_i(x^i,\theta^i)={\theta^i}^\top \varphi^n(x^i)-\lambda\|\theta^i\|^2$
is $2\lambda$-strongly concave in $\theta^i$ on the closed convex set $\hat{\Theta}_N$. Hence, the lower-level problem admits a unique solution
$\theta^{*}(x^i)
 =
\mathrm{Proj}_{\hat{\Theta}_N}
\Bigl(\tfrac{1}{2\lambda}\varphi^n(x^i)\Bigr)$ for $i=1,2$,
where $\mathrm{Proj}_S(\cdot)$ denotes the projection onto the set $S$. Moreover,
\[
u_i(x^i,\theta^{*}(x^i))-u_i(x^i,\theta)
\geq
\lambda\|\theta-\theta^{*}(x^i)\|^2
\geq
\lambda\|\theta-\theta^{*}(x^i)\|_\infty^2.
\]
Therefore, the lower-level growth condition introduced in Section 4 holds with $p_\lambda=2$ and $\gamma_\lambda=\lambda$, which explains the stabilizing effect of regularization on the lower-level worst-case parameter selection.
}
Substituting $\theta^{*}(x^i)$ into the upper-level problem simplifies the problem \eqref{num-eq-1} to:
\begin{equation}
0 \in F(x)
+  \mathcal{N}_{\mathcal{X}_1\times\mathcal{X}_2}(x), \label{num-final}
\end{equation}
where
\[
F(x) =
\begin{pmatrix}
F_1(x)\\[4pt]
F_2(x)
\end{pmatrix},
\quad
F_i(x)
= \sum_{j=1}^n \theta^{*}_j(x^i)\mathbb{E}_{Q_j}[\nabla_{x^i}\varphi(x^i,\xi)]
  +\kappa\bigl(x^i-x^{-i}\bigr),\
 i=1,2.
\]
For problem \eqref{num-final}, we have the following conclusion:
\begin{proposition}\label{prop:danskin}
Assume that for each $j=1,\ldots,n$, (i) there exists $\beta\in[0,1)$ such that for each $i=1,2$ and all $x^i\in\mathcal{X}_i$, ${x^i}^\top\xi\ge-\beta$ almost surely under \(Q_j\), and (ii) there exists $\bar{\sigma}_j^2<\infty$ such that $\mathbb{E}_{Q_j}\|\xi\|^2\leq\bar{\sigma}_j^2$.
      Then, $F$
      is monotone and Lipschitz continuous on $\mathcal X_1\times\mathcal{X}_2$ with modulus
      \(
        L_F
        =(\sqrt{2}+\frac{\sqrt{2}n}{2\lambda})\frac{\max_j\bar{\sigma}_j^2}{(1-\beta)^2}+2\sqrt{2}\kappa.
      \)
\end{proposition}

 \textit{Proof}.
For each $x\in\mathcal X_1\times\mathcal{X}_2$ and direction $d=(d^1,d^2)\in\mathbb R^{2d}$,
define
\( q(t):=\bigl[f_\lambda(x+td)-f_\lambda(x)\bigr]/t
\), we have
\begin{equation*}
    \begin{aligned}
        \limsup_{t\downarrow0}q(t)&=\limsup_{t\downarrow0}\frac{f_\lambda(x+td)-f_\lambda(x)}{t}\\
        &\geq\sum_{i=1}^2\limsup_{t\downarrow0}\frac{U_i(x^i+td^i)-U_i(x^i)}{t}+\kappa(x^1-x^2)^\top(d^1-d^2)\\
        &\geq\sum_{i=1}^2\limsup_{t\downarrow0}\frac{u_i(x^i+td^i,\theta^{*}(x^i))-u_i(x^i,\theta^{*}(x^i))}{t}+\kappa(x^1-x^2)^\top(d^1-d^2)\\
        &=\sum_{i=1}^2\nabla_{x^i}u_i(x^i,\theta^{*}(x^i))^\top d^i+\kappa(x^1-x^2)^\top(d^1-d^2).\\
    \end{aligned}
\end{equation*}
Similarly, we have
$
        \liminf_{t\downarrow0}q(t)\leq
\sum_{i=1}^2\nabla_{x^i}u_i(x^i,\theta^{*}(x^i))^\top d^i+\kappa(x^1-x^2)^\top(d^1-d^2).$
Thus, the directional derivative exists. Moreover,
by Danskin's theorem, we have
$$
\nabla_{x^i} f_\lambda(x)=\nabla_{x^i}u_i(x^i,\theta^{*}(x^i))+\kappa(x^i-x^{-i})        =\sum_{j=1}^{n}\theta^{*}_j(x^i) \mathbb E_{Q_j}[\nabla_{x^i}\varphi(x^i,\xi)]+\kappa(x^i-x^{-i})  =F_i(x).
$$
With the convexity of $f_\lambda$, we can deduce that $F=\nabla f_\lambda$ is monotone.

For the Lipschitz property, denote $g_j(x^i)=\mathbb{E}_{Q_j}[\nabla_{x^i}\varphi(x^i,\xi)]$, we have from the definition of $\varphi$ that
\[
\|g_j(x^i)\|\leq\frac{\bar{\sigma}_j}{(1-\beta)}=:L_{\varphi,j}, \ \|\nabla g_j(x^i)\|\leq\mathbb{E}_{Q_j}\left[\frac{\|\xi\|^2}{(1-\beta)^2}\right]\leq\frac{\bar{\sigma}_j^2}{(1-\beta)^2}=:L_{g,j}.
\]
Thus $g_j$ is $L_{g,j}$-Lipschitz, $\mathbb E_{Q_j}[\varphi(\cdot,\xi)]$ is $L_{\varphi,j}$-Lipschitz and then $\theta^{*}(\cdot)$ is $\frac{\sqrt{n}\max_jL_{\varphi,j}}{2\lambda}$-Lipschitz.
Based on the above results,
we have
\begin{equation*}
    \begin{aligned}
      \|F_i(x)-F_i(y)\|
 &= \|\sum_j \theta^{*}_j(x^i)g_j(x^i)-\theta^{*}_j(y^i)g_j(y^i)+\kappa(x^i-y^i)-\kappa(x^{-i}-y^{-i})\|\\
 &\leq\sum_j \theta^{*}_j(x^i)\|g_j(x^i)-g_j(y^i)\|+\sum_j \|g_j(y^i)\| |\theta^{*}_j(x^i)-\theta^{*}_j(y^i)|+\kappa\|x^{-i}-y^{-i}\|+\kappa\|x^{i}-y^{i}\|\\
 &\leq \max_jL_{g,j}\|x^i-y^i\|+\frac{n\max_jL_{\varphi,j}^2}{2\lambda}\|x^i-y^i\|+2\kappa\|x-y\|\leq L_{\bar F}\|x-y\|,
\end{aligned}
\end{equation*}
where $L_{\bar F}:=(1+\frac{n}{2\lambda})\frac{\max_j\bar{\sigma}_j^2}{(1-\beta)^2}+2\kappa$.
By setting $L_F=\sqrt{2}L_{\bar F}$, the stated modulus is proved.
$\hfill\blacksquare$

As for the assumptions in Proposition \ref{prop:danskin}, assumption (i) is mild in financial contexts where  ${x^i}^\top\xi\geq-1$ is often satisfied for $i=1,2$.
{\color{black}
For any fixed $\lambda>0$, Proposition~\ref{prop:danskin} shows that the operator $F$ in \eqref{num-final} is monotone and $L_F$-Lipschitz continuous on the nonempty closed convex set $\mathcal{X}_1\times\mathcal{X}_2$. Therefore, the standard convergence theory of the extragradient method for monotone Lipschitz variational inequalities applies to the fixed-$\lambda$ problem \eqref{num-final}; see, e.g., \cite{nemirovski2004prox,malitsky2015projected}.
In particular, Algorithm~\ref{alg:extragradient} is covered by this theory whenever the constant step size is chosen such that $\eta\in\bigl(0,\frac{1}{2L_F}\bigr)$. Moreover,
define the variational inequality
gap function by
$\operatorname{Gap}(z):=\max_{u\in \mathcal{X}_1\times\mathcal{X}_2}
\langle F(u),\, z-u\rangle .$
Then, from \cite{juditsky2011solving,nemirovski2004prox}, with a constant stepsize
$\eta$ of order $1/L_F$, the averaged auxiliary iterates $\bar y^K:=\frac1K\sum_{k=0}^{K-1} y^k$
satisfy
\[
\operatorname{Gap}(\bar y^K)\le \frac{C\,L_F}{K},
\]
where $C>0$ is a constant independent of $K$ and $L_F$. Therefore, to achieve
$\operatorname{Gap}(\bar y^K)\le \varepsilon$, it suffices to take
$K=O\!\left(\frac{L_F}{\varepsilon}\right).$
Since Proposition~\ref{prop:danskin} yields $L_F=O(1/\lambda)$ as $\lambda\downarrow0$,
the required number of extragradient iterations scales as
$K=O\!\left(\frac{1}{\lambda\varepsilon}\right).$
}

\begin{algorithm}
\caption{Extragradient method for the fixed-$\lambda$ regularized Bayesian DRVI problem \eqref{num-final}}
\label{alg:extragradient}
\begin{algorithmic}[1] 
\Require {\color{black}Regularization parameter $\lambda>0$}, step size $\eta \in \bigl(0, \frac{1}{2L_F} \bigr)$, tolerance $\varepsilon > 0$, initial point $x^{0} \in \mathcal{X}_1\times\mathcal{X}_2$
\For{$k = 0, 1, 2, \dots$}
    \State \textbf{Compute $F$ at $x^k$:} \For{$i=1,2$}
        \State $z^{i,k} \gets \frac{1}{2\lambda}\varphi^n(x^{i,k})$
        \State $\theta^{\lambda,i} \gets \operatorname{Proj}_{\hat{\Theta}_N}(z^{i,k})$
        \State $F_i^k \gets \sum_{j=1}^n \theta_j^{\lambda,i} \mathbb{E}_{Q_j}[\nabla_{x^i}\varphi(x^{i,k},\xi)] + \kappa(x^{i,k}-x^{-i,k})$
    \EndFor
    \State $F^k \gets (F_1^k, F_2^k)^\top$
 \State \textbf{Extragradient step:}
    $y^k \gets \operatorname{Proj}_{\mathcal{X}_1\times\mathcal{X}_2}\bigl(x^k - \eta F^k\bigr)$
    \State \textbf{Recompute $F$ at $y^k$:}
    \For{$i=1,2$}
        \State $\tilde{z}^{i,k} \gets \frac{1}{2\lambda}\varphi^n(y^{i,k})$
        \State $\tilde{\theta}^{\lambda,i} \gets \operatorname{Proj}_{\hat{\Theta}_N}(\tilde{z}^{i,k})$
        \State $\tilde{F}_i^k \gets \sum_{j=1}^n \tilde{\theta}_j^{\lambda,i} \mathbb{E}_{Q_j}[\nabla_{x^i}\varphi(y^{i,k},\xi)] + \kappa(y^{i,k}-y^{-i,k})$
    \EndFor
    \State $\tilde{F}^k \gets (\tilde{F}_1^k, \tilde{F}_2^k)^\top$
      \State \textbf{Update:} $x^{k+1} \gets \operatorname{Proj}_{\mathcal{X}_1\times\mathcal{X}_2}\bigl(x^k - \eta \tilde{F}^k\bigr)$ and $  \bar y^{\,k+1}\gets \frac{1}{k+1}\sum_{t=0}^{k}y^t$
    \If{$\operatorname{Gap}(\bar y^{\,k+1})\le \varepsilon$}
        \State \Return $x^*\gets \bar y^{\,k+1}$
    \EndIf
\EndFor
\end{algorithmic}
\end{algorithm}

To further specify the setting of problem \eqref{num-eq-1}, we consider $d=10$ risky assets with $n=3$, representing three potential market regimes: bear, oscillating and bull markets. To generate the necessary data for stock returns, we
adopt the following factor model widely used in \cite{chen2018data,zhu2014portfolio}:
\[\xi_m^j=a_m^j+b_m^j\zeta^j_{market}+e^j_m,\ j=1,2,3,\ m=1,\ldots,10,\]
where the scalar market factors follow
\(\zeta_{\mathrm{market}}^j\sim\mathcal{N}(\mu_j,\sigma_{z,j}^2)\) , $j=1,2,3$, with
\[
(\mu_1,\sigma_{z,1}^2)=(0.005,0.001),\quad
(\mu_2,\sigma_{z,2}^2)=(0.010,0.004),\quad
(\mu_3,\sigma_{z,3}^2)=(0.015,0.007).
\]
The intercepts \(a_m^j\) capture each asset’s baseline return and are set to \(-0.0005\) for bear market $(j=1)$, \(0\) for oscillating market $(j=2)$, and \(0.0005\) for bull market $(j=3)$; the factor loadings \(b_m^j\) determine sensitivity to the market factor and increase linearly with \(m\):
\(b_m^1=0.60+0.02(m-1)\) in bear market,
\(b_m^2=0.80+0.015(m-1)\) in oscillating market, and
\(b_m^3=1.00+0.01(m-1)\) in bull market. The residuals
\(e_m^j\sim\mathcal{N}(0,\sigma_{e,j}^2)\) follow independent normal distributions with variances \(\sigma_{e,1}^2=2\times10^{-5}\),
\(\sigma_{e,2}^2=1\times10^{-5}\), and
\(\sigma_{e,3}^2=0.5\times10^{-5}\), reflecting higher asset-specific risk in downturns and lower risk in upturns. The mixture weight is set to \(\theta^c=(0.5,0.25,0.25)\). Here, we set $\kappa=0.1$, $\lambda=0.01$ and \( r^c = 0\) to model a risk-neutral decision-maker.

We generate a large dataset $\{\xi^s\}_{s=1}^{N_{\rm sim}}$ with $N_{\rm sim}=10^6$ to simulate the true equilibrium solution $x^c$ as a benchmark. For comparison, we compute the approximated solution $x^*$ by different data-driven methods on training samples with size $N\in\{ 20,50,500,1000,3000\}$:
{\color{black}
(i) the empirical SVI solved by sample average approximation (SAA) \cite{xu2010sample};
(ii) the posterior predictive SVI, which replaces the true distribution by the Bayesian posterior predictive distribution \cite{shapiro2025episodic,wu2018bayesian} serving as a non-robust Bayesian benchmark;
(iii) two frequentist robust benchmarks with $\ell_1$ and modified $\chi^2$ ambiguity sets; and
(iv) our Bayesian DRVI method.
For the $\ell_1$-, modified $\chi^2$-, and our Bayesian DRVI methods, the ambiguity sets are calibrated under the same level $\alpha=0.05$.
}
Each method is run over 100 independent trials. For each trial, to quantify the accuracy of the solution $x^*$,  we compute the residual distance: $d(x^{*})=\|x^{*}-x^{c}\|$ and the squared relative utility error \cite{wu2018bayesian}:
$$e(x^*)=\left(\frac{U(x^*)-U(x^c)}{U(x^c)}\right)^2,\text{ where } U(x)= \frac1{N_{\rm sim}}\sum_{s=1}^{N_{\rm sim}} -\log\bigl(1+{\xi^s}^\top x^{1}\bigr)+\frac{\kappa}{2}\|x^{1}-x^{2}\|^2.$$
Here, since the game is symmetric across the two accounts, we focus on account 1’s utility without loss of generality.
Furthermore, to quantify out-of-sample performance, we generate a fixed test set $\{\xi^t\}_{t=1}^{N_{\rm test}}$ with \(N_{\rm test}=5000\). Based on the single-period returns \(\{R_t = {\xi^t}^\top x^{1*}\}_{t=1}^{N_{\rm test}}\), we further assess the out-of-sample performance of the optimal portfolio with additional metrics in  \cite{bielecki2018unified}--Sharpe ratio (measuring return per unit of total risk) and risk-adjusted return on capital (RAROC; measuring return per unit of tail risk) at $10\%$ tail level.
Specifically, let the out-of-sample mean and standard deviation of the returns be
\[
\bar R=\frac{1}{N_{\rm test}}\sum_{t=1}^{N_{\rm test}} R_t,\qquad
\hat\sigma_R=\sqrt{\frac{1}{N_{\rm test}-1}\sum_{t=1}^{N_{\rm test}}\bigl(R_t-\bar R\bigr)^2}.
\]
The Sharpe ratio is then
\[
\widehat{\mathrm{SR}}=\frac{\bar R}{\hat\sigma_R}.
\]
For the RAROC at $10\%$ tail level, define the losses $l_t:=-R_t$.
Let the empirical value-at-risk be $\widehat{\mathrm{VaR}}_{10\%}$ as the $10\%$-quantile of $\{l_t\}_{t=1}^{N_{\rm test}}$, and let the empirical  expected shortfall of losses be
\[
\widehat{\mathrm{ES}}_{10\%}=\frac{1}{\#\{t:l_t\ge \widehat{\mathrm{VaR}}_{10\%}\}}\sum_{t:l_t\ge \widehat{\mathrm{VaR}}_{10\%}} l_t,
\]
where $\#\{\cdot\}$ denotes the cardinality (number of elements) of a set.
Then the RAROC at $10\%$ tail level is calculated as
\[
\widehat{\mathrm{R}}_{10\%}=\frac{\bar R}{\widehat{\mathrm{ES}}_{10\%}}.
\]

We report in Tables \ref{tab:residuals_scaled}-\ref{tab:raroc_scaled} the sample means and variances of all performance metrics across 100 trials, allowing a comprehensive assessment of each method’s estimation accuracy, sensitivity to sampling variability, and risk-adjusted performance as the training sample size increases.

\begin{table}
  \caption{{\color{black}Residual distance $d(x^*)$}: means and variances across training sample sizes\label{tab:residuals_scaled}}{%
    \setlength{\tabcolsep}{7pt}
    \begin{tabular}{l l c c c c c}
      \hline
      & Method & \(N=20\) & \(N=50\) & \(N=500\) & \(N=1000\) & \(N=3000\) \\
      \hline
      \multirow{5}{*}{Mean (\(\times10^{-1}\))} & \(\ell_1\) norm        & 7.311 & 6.914 & 5.629 & 3.203 & 0.669 \\
                                                & Modified \(\chi^2\)    & 7.111 & 7.027 & 6.456 & 6.336 & 6.221 \\
                                                & SAA                    & 6.205 & 5.325 & 1.758 & 0.805 & 0.217 \\
                                                & Bayesian               & \textbf{6.123} & \textbf{5.019} & \textbf{1.555} & \textbf{0.775} & \textbf{0.158} \\
                                                & \textcolor{black}{Posterior predictive} & \textcolor{black}{6.182} & \textcolor{black}{5.266} & \textcolor{black}{1.694} & \textcolor{black}{0.806} & \textcolor{black}{0.214} \\
      \hline
      \multirow{5}{*}{Variance (\(\times10^{-2}\))} & \(\ell_1\) norm        & 4.2579 & 4.1647 & 3.4892 & 3.2749 & 2.5608 \\
                                                    & Modified \(\chi^2\)    & \textbf{3.6275} & \textbf{3.2301} & \textbf{2.2205} & \textbf{1.2423} & 0.7599 \\
                                                    & SAA                    & 4.9777 & 4.2403 & 3.5268 & 2.4571 & 0.8582 \\
                                                    & Bayesian               & 4.1647 & 4.0356 & 3.3361 & 2.2159 & \textbf{0.7046} \\
                                                    & \textcolor{black}{Posterior predictive} & \textcolor{black}{4.7948} & \textcolor{black}{4.1757} & \textcolor{black}{3.5347} & \textcolor{black}{2.2902} & \textcolor{black}{0.8152} \\
      \hline
    \end{tabular}%
  }{}
\end{table}

\begin{table}
  \caption{{\color{black}Squared relative utility error $e(x^*)$}: means and variances across training sample sizes\label{tab:util_err_scaled}}{%
    \setlength{\tabcolsep}{7pt}
    \begin{tabular}{l l c c c c c}
      \hline
      & Method & \(N=20\) & \(N=50\) & \(N=500\) & \(N=1000\) & \(N=3000\) \\
      \hline
      \multirow{5}{*}{Mean (\(\times10^{-3}\))} & \(\ell_1\) norm        & 3.843 & 3.440 & 2.625 & 1.636 & 0.368 \\
                                                & Modified \(\chi^2\)    & 3.300 & 3.973 & 4.024 & 4.095 & 4.273 \\
                                                & SAA                    & 1.464 & 1.210 & 1.034 & 0.330 & 0.151 \\
                                                & Bayesian               & \textbf{1.360} & \textbf{1.113} & \textbf{0.828} & \textbf{0.315} & \textbf{0.125} \\
                                                & \textcolor{black}{Posterior predictive} & \textcolor{black}{1.423} & \textcolor{black}{1.181} & \textcolor{black}{1.044} & \textcolor{black}{0.330} & \textcolor{black}{0.151} \\
      \hline
      \multirow{5}{*}{Variance (\(\times10^{-6}\))} & \(\ell_1\) norm        & 2.3010 & 1.9958 & 1.5284 & 1.4141 & 0.8297 \\
                                                    & Modified \(\chi^2\)    & 2.3468 & \textbf{1.0422} & \textbf{0.5144} & \textbf{0.2289} & 0.1330 \\
                                                    & SAA                    & 1.8827 & 1.2819 & 0.8018 & 0.4076 & 0.0649 \\
                                                    & Bayesian               & \textbf{1.5179} & 1.2662 & 0.7540 & 0.3783 & \textbf{0.0367} \\
                                                    & \textcolor{black}{Posterior predictive} & \textcolor{black}{1.8601} & \textcolor{black}{1.2795} & \textcolor{black}{0.7946} & \textcolor{black}{0.4071} & \textcolor{black}{0.0653} \\
      \hline
    \end{tabular}%
  }{}
\end{table}

\begin{table}
  \caption{{\color{black}Sharpe ratio $\widehat{\mathrm{SR}}$}: means and variances across training sample sizes\label{tab:sharpe_scaled}}{%
    \setlength{\tabcolsep}{7pt}
    \begin{tabular}{l l c c c c c}
      \hline
      & Method & \(N=20\) & \(N=50\) & \(N=500\) & \(N=1000\) & \(N=3000\) \\
      \hline
      \multirow{5}{*}{Mean (\(\times10^{-1}\))} & \(\ell_1\) norm        & 1.798 & 1.794 & 1.800 & 1.808 & 1.821 \\
                                                & Modified \(\chi^2\)    & 1.795 & 1.791 & 1.789 & 1.790 & 1.789 \\
                                                & SAA                    & 1.826 & \textbf{1.828} & \textbf{1.829} & \textbf{1.827} & \textbf{1.826} \\
                                                & Bayesian               & 1.820 & 1.824 & \textbf{1.829} & \textbf{1.827} & \textbf{1.826} \\
                                                & \textcolor{black}{Posterior predictive} & \textcolor{black}{\textbf{1.827}} & \textcolor{black}{\textbf{1.828}} & \textcolor{black}{\textbf{1.829}} & \textcolor{black}{\textbf{1.827}} & \textcolor{black}{\textbf{1.826}} \\
      \hline
      \multirow{5}{*}{Variance (\(\times10^{-6}\))} & \(\ell_1\) norm        & 1.7977 & 1.6436 & 0.8393 & 1.1894 & 0.7318 \\
                                                    & Modified \(\chi^2\)    & \textbf{1.2556} & \textbf{0.4248} & \textbf{0.1006} & \textbf{0.0889} & \textbf{0.0738} \\
                                                    & SAA                    & 2.0414 & 1.0569 & 0.2636 & 0.1862 & 0.1605 \\
                                                    & Bayesian               & 1.7583 & 0.9653 & 0.2527 & 0.1542 & 0.1268 \\
                                                    & \textcolor{black}{Posterior predictive} & \textcolor{black}{1.9704} & \textcolor{black}{1.0689} & \textcolor{black}{0.2611} & \textcolor{black}{0.1941} & \textcolor{black}{0.1607} \\
      \hline
    \end{tabular}%
  }{}
\end{table}

\begin{table}
  \caption{{\color{black}RAROC at 10\% tail level $\widehat{\mathrm{R}}_{10\%}$}: means and variances across training sample sizes\label{tab:raroc_scaled}}{%
    \setlength{\tabcolsep}{7pt}
    \begin{tabular}{l l c c c c c}
      \hline
      & Method & \(N=20\) & \(N=50\) & \(N=500\) & \(N=1000\) & \(N=3000\) \\
      \hline
      \multirow{5}{*}{Mean (\(\times10^{-1}\))} & \(\ell_1\) norm        & 1.163 & 1.160 & 1.164 & 1.171 & 1.180 \\
                                                & Modified \(\chi^2\)    & 1.160 & 1.157 & 1.155 & 1.155 & 1.154 \\
                                                & SAA                    & \textbf{1.184} & \textbf{1.185} & \textbf{1.185} & \textbf{1.184} & \textbf{1.183} \\
                                                & Bayesian               & 1.180 & 1.183 & \textbf{1.185} & \textbf{1.184} & \textbf{1.183} \\
                                                & \textcolor{black}{Posterior predictive} & \textcolor{black}{\textbf{1.184}} & \textcolor{black}{\textbf{1.185}} & \textcolor{black}{\textbf{1.185}} & \textcolor{black}{\textbf{1.184}} & \textcolor{black}{\textbf{1.183}} \\
      \hline
      \multirow{5}{*}{Variance (\(\times10^{-6}\))} & \(\ell_1\) norm        & 1.2628 & 1.1974 & 0.6883 & 0.5591 & 0.1941 \\
                                                    & Modified \(\chi^2\)    & \textbf{0.9043} & \textbf{0.3301} & \textbf{0.0989} & \textbf{0.0881} & \textbf{0.0623} \\
                                                    & SAA                    & 1.6741 & 0.9903 & 0.1419 & 0.0941 & 0.0675 \\
                                                    & Bayesian               & 1.2089 & 0.6837 & 0.1366 & 0.0908 & 0.0629 \\
                                                    & \textcolor{black}{Posterior predictive} & \textcolor{black}{1.3723} & \textcolor{black}{0.8958} & \textcolor{black}{0.1401} & \textcolor{black}{0.0936} & \textcolor{black}{0.0668} \\
      \hline
    \end{tabular}%
  }{}
\end{table}

The observations derived from Tables \ref{tab:residuals_scaled}-\ref{tab:raroc_scaled} lead to the following conclusions.

\begin{itemize}
\item \textbf{Consistency and estimation accuracy.}
Across all sample sizes, the Bayesian DRVI attains the smallest mean residual distances to the benchmark $x^c$ and the lowest squared relative utility errors (Tables \ref{tab:residuals_scaled}-\ref{tab:util_err_scaled}), demonstrating its superior estimation accuracy. Moreover, Table \ref{tab:residuals_scaled} shows that the Bayesian solutions improve steadily as $N$ increases, consistent with the theoretical convergence results. In contrast, the $\ell_1$ norm method is outperformed by Bayesian DRVI, posterior predictive SVI, and empirical SVI, while the modified $\chi^2$ scheme exhibits a persistent bias, indicative of over-conservatism that does not vanish with additional data.

\item \textbf{Stability and variance reduction.}
Tables \ref{tab:residuals_scaled}-\ref{tab:util_err_scaled} show that the Bayesian DRVI maintains low variance across sample sizes, particularly for large $N$, and achieves the lowest variance among all methods for $N=3000$. {\color{black}
The posterior predictive SVI is also more stable than SAA in most cases, which is consistent with the fact that it remains within a fixed parametric model rather than relying on the empirical distribution. However, because it incorporates posterior information only through Bayesian averaging under a single predictive distribution, it is still less stable than Bayesian DRVI, which explicitly accounts for such uncertainty through posterior-credible robustification.
}
While the modified $\chi^2$ approach also displays small variances, this is accompanied by pronounced bias, reflecting variance reduction driven by excessive pessimism rather than genuine robustness. In contrast, among the essentially unbiased approaches ($\ell_1$ norm, Bayesian and SAA), the Bayesian method consistently attains both the lowest mean errors and the smallest variances, highlighting its superior stability and reliability in finite samples.

\item \textbf{Risk-return performance.}
The SAA and Bayesian methods form a top tier of performance for $N \geq 500$, achieving nearly identical and superior results in both the Sharpe ratio and RAROC at $10\%$ tail level (Tables \ref{tab:sharpe_scaled}-\ref{tab:raroc_scaled}). In smaller samples ($N=20,50$), SAA may achieve a slightly higher mean performance, but the Bayesian method typically exhibits smaller variance, confirming that it provides more stable risk-return trade-offs across different data realizations.
The inherent conservatism of the $\ell_1$ and modified $\chi^2$ methods leads to suboptimal decisions, resulting in consistently lower performance.
This demonstrates the substantial opportunity cost of using overly pessimistic ambiguity sets, while the Bayesian approach delivers excellent risk-adjusted and tail-controlled returns with improved reliability.
\end{itemize}

In summary, the proposed Bayesian DRVI method effectively balances the trade-off between robustness and optimality.
{\color{black}
Compared with posterior predictive and empirical SVI models, the proposed method demonstrates the additional value of posterior robustification beyond Bayesian updating alone and provides more reliable decisions. Compared with the frequentist robust benchmarks, it also mitigates excessive conservatism.
}

\section{Concluding remarks.}

Building on the DRVI framework in \cite{sun2023distributionally}, we consider a Bayesian formulation \eqref{bdrvi-true} which captures distributional uncertainty through the parametric ambiguity set \(\Theta^c\).
To ensure tractability, we construct a data-driven ambiguity set with posterior coverage guarantees and introduce the regularized Bayesian DRVI model \eqref{bdrvi-regular2}. Under mild assumptions, we prove the existence of solutions for both Bayesian DRVI \eqref{bdrvi2} and its regularized counterpart \eqref{bdrvi-regular2}.
We further establish the convergence of the solutions of \eqref{bdrvi-regular2} to a solution of the problem \eqref{bdrvi-true} as the sample size grows and the regularization vanishes.
Moreover, as \(\lambda\downarrow 0\), the regularized maximizer of the lower-level problem in \eqref{bdrvi-regular2} converges to the least-norm maximizer of the lower-level problem in \eqref{bdrvi2}.
To address practical finite-sample and contamination effects, we derive finite-sample guarantees for the optimal solutions of \eqref{bdrvi-regular2} and quantify the solution stability between this formulation and its perturbed variant \eqref{p-bdrvi2}.
The effectiveness of the proposed framework is demonstrated through its application to a distributionally robust multi-portfolio Nash equilibrium problem. By integrating Bayesian learning with DRVI, our framework provides a robust and practical foundation for decision-making under uncertainty.

It is natural to extend the Bayesian DRVI framework to stochastic generalized equations. While statistical robustness results for stochastic generalized equations have been established in \cite{guo2023data}, analogous theoretical guarantees for the Bayesian DRVI framework remain an open question. We leave these extensions to future work.

\section*{Acknowledgments.}
This work was supported in part by grants from Hong Kong Research Grants Council (PolyU15300123, PolyU15300124) and National Key R\&D Program of China (2022YFA1004000, 2022YFA1004001).

\bibliographystyle{siam}
\bibliography{ref}
\end{document}